\documentclass[a4paper,12pt,reqno]{amsart}
\usepackage{amsmath}
\usepackage{amssymb}
\usepackage{amsthm}

\usepackage{amscd}
\usepackage{amsfonts}
\usepackage{setspace}
\usepackage{version}

\usepackage{hyperref}

\title[Barban--Davenport--Halberstam asymptotics for general sequences]{Simple Barban--Davenport--Halberstam type asymptotics for general sequences}
\author{Adam J Harper}
\address{Mathematics Institute, Zeeman Building, University of Warwick, Coventry CV4 7AL, England}
\email{A.Harper@warwick.ac.uk}
\date{27th December 2024}
\thanks{Quite a lot of the research leading to this paper (in the particular case of the smooth numbers) was carried out around 2013, when the author was supported by a postdoctoral fellowship from the Centre de Recherches Math\'ematiques, Montr\'eal. It was completed more recently, with partial support from the Engineering and Physical Sciences Research Council of the United Kingdom [grant EP/V055755/1]; and with support from the Swedish Research Council [grant no. 2021-06594] while the author was in residence at Institut Mittag-Leffler, Djursholm, Sweden, during the 2024 Analytic Number Theory program.}

\numberwithin{equation}{section}
\theoremstyle{plain}

\onehalfspace
\newcommand{\N}{\mathbb{N}}

\newtheorem{thm1}{Theorem}
\newtheorem{thm2}[thm1]{Theorem}
\newtheorem{cor1}{Corollary}

\newtheorem{lem1}{Lemma}
\newtheorem{lem2}[lem1]{Lemma}
\newtheorem{lem3}[lem1]{Lemma}

\begin{document}

\maketitle

\begin{abstract}
We prove two estimates for the Barban--Davenport--Halberstam type variance of a general complex sequence in arithmetic progressions. The proofs are elementary, and our estimates are capable of yielding an asymptotic for the variance when the sequence is sufficiently nice, and is either somewhat sparse or is sufficiently like the integers in its divisibility by small moduli.

As a concrete application, we deduce a Barban--Davenport--Halberstam type variance asymptotic for the $y$-smooth numbers less than $x$, on a wide range of the parameters. This addresses a question considered by Granville and Vaughan.
\end{abstract}

% INTRODUCTION %%%%%%%%%%%%%%%%%%%%%%%%%%%%%
\section{Introduction}
In this note we are interested in the variance of the distribution of sets and sequences in arithmetic progressions. More specifically, suppose that $x$ is large and that $\mathcal{A} = (a_n)_{n \leq x}$ is a complex sequence, and for $z \leq x$ consider the following counting functions:
$$ \mathcal{A}(z) := \sum_{n \leq z} a_n , \;\;\;\;\; \mathcal{A}(z;q,a) := \sum_{n \leq z, n \equiv a \; \textrm{mod } q} a_n . $$
We will also need notation for the average of $\mathcal{A}(z;q,a)$ over values $a$ such that $(a,q)=h$:
$$ \text{Aver}(\mathcal{A},z;q,h) := \frac{1}{\phi(q/h)} \sum_{1 \leq a \leq q, (a,q)=h} \mathcal{A}(z;q,a), \;\;\; h \mid q. $$
Our goal will be to give an asymptotic formula (i.e. one or more main terms, and some smaller order error terms) for the following quantity:
$$ V(\mathcal{A},x,Q) := \sum_{Q/2 < q \leq Q} \sum_{1 \leq a \leq q} \left|\mathcal{A}(x;q,a) - \text{Aver}(\mathcal{A},x;q,(a,q)) \right|^{2}. $$
In particular, we aim to obtain an asymptotic on a wide range of $Q$, and making the weakest and simplest assumptions that we can on $\mathcal{A}$.

\vspace{12pt}
Before stating our results we briefly describe some of the extensive history of this problem. Probably the most studied specific case is where $\mathcal{A}$ is (the indicator function of) the set of primes, or the values of the von Mangoldt function $\Lambda(n)$. Barban (see \cite{barban}) and Davenport and Halberstam~\cite{davhal} proved upper bounds for (essentially) $V(\Lambda,x,Q)$, which were sharpened by Gallagher~\cite{gallagher} to yield the following result: for any fixed $A > 0$, and any $x \log^{-A}x \leq Q \leq x$, we have
$$ \sum_{q \leq Q} \sum_{(a,q)=1} \left| \Psi(x;q,a)-\frac{x}{\phi(q)} \right|^{2} \ll_{A} xQ \log x , $$
where $\Psi(z;q,a) := \sum_{n \leq z, n \equiv a (\textrm{mod } q)} \Lambda(n)$. Note that if $(a,q) \neq 1$ then $\Psi(x;q,a) \ll \log x$ (and typically $\Psi(x;q,a) =0 $), so such values would make a negligible contribution.

Gallagher's bound implies that $| \Psi(x;q,a)-\frac{x}{\phi(q)} |$ is $O(\sqrt{(x\log x)/\phi(q)})$ on average, which is the order one would expect if primes were distributed amongst the coprime residue classes ``at random''. Thus one might expect that Gallagher's upper bound is asymptotically correct. This was confirmed by Montgomery~\cite{mont2}, and then Hooley~\cite{hooley} proved the following sharper result: for any fixed $A > 0$, and any $Q \leq x$,
$$ \sum_{q \leq Q} \sum_{(a,q)=1} \left| \Psi(x;q,a)-\frac{x}{\phi(q)} \right|^{2} = Qx\log Q + B_{1}Qx + O(Q^{5/4}x^{3/4}) + O_{A}(x^{2}/\log^{A}x) , $$
where $B_{1}$ is an explicit constant. Several methods are used in the works mentioned, including the multiplicative large sieve, average results in additive (Goldbach-type) prime number theory, and an important divisor switching idea that we shall explain later. One also uses Siegel--Walfisz type information about primes in progressions to small moduli, which produces the error term $O_{A}(x^{2}/\log^{A}x)$ and means that one only gets a genuine asymptotic for the variance when $Q \geq x \log^{-A}x$.

Goldston and Vaughan~\cite{goldstonvaughan} used an exponential sums argument to improve the first ``big Oh'' term in the above display to $O(Q^{3/2}x^{1/2})$. More precise theorems can be proved assuming the Generalised Riemann Hypothesis, and much wider ranges of $Q$ can be handled if one seeks lower bounds (rather than upper bounds or asymptotics) for the variance. One can also prove results involving higher moments (i.e. higher powers of the discrepancy); where the variables $q,a$ run over restricted sets; etc.. See e.g. Hooley's survey article~\cite{hooleysurv} for further details of much of the work in the 20th century, and the exciting recent papers of La Bret\`eche and Fiorilli~\cite{dlbfiormom} and of Fiorilli and Martin~\cite{fiormartdis} (and references there) for subsequent developments. 

\vspace{12pt}
Quite a lot has also been done to obtain Barban--Davenport--Halberstam type bounds, and then asymptotics, in an ``axiomatic'' way for general $\mathcal{A}$. The state of the art is perhaps due to Vaughan~\cite{vaughangbdh}, (see also Hooley's later paper~\cite{hooley14}), who assumes that
$$ \mathcal{A}(z;q,a) = z f(q,(a,q)) + O\left(\frac{z}{\Psi(z)}\right) \;\;\; \forall z \geq 1,\;\;\; \forall q, a \in \N, $$
for some increasing function $\Psi(z)$ such that $\Psi(z) > \log z$ and $\int_{1}^{z} \Psi(y)^{-1} dy \ll z/\Psi(z)$. This generalises the so-called ``Criterion U'' of Hooley~\cite{hooley9}, in which the ``big Oh'' term is assumed to be $O_{A}(z/\log^{A}z)$ for every fixed $A > 0$. Vaughan also assumes a uniform bound $\sum_{n \leq z} |a_n|^2 \ll z$ for all $z \geq 1$. Under these conditions, Vaughan~\cite{vaughangbdh} used an exponential sum method (developing that of Goldston and Vaughan~\cite{goldstonvaughan}) to show that if $x \geq Q > \sqrt{x} \log 2x$, then
\begin{eqnarray}
 \sum_{q \leq Q} \sum_{1 \leq a \leq q} \left|\mathcal{A}(x;q,a) - x f(q,(a,q)) \right|^{2} & = & Q\sum_{n \leq x} |a_n|^2 - Qx \sum_{q=1}^{\infty} g(q) + O\left(Q^{2} \int_{1}^{x/Q} (\sum_{q > y} g(q)) dy \right) \nonumber \\
&& + O\left(x^{3/2}\log x + \frac{x^{2}(\log 2x)^{9/2}}{\Psi(x)} + \frac{x^{2}(\log x)^{4/3}}{\Psi(x)^{2/3}} \right) . \nonumber
\end{eqnarray}
Here $g(q) := \phi(q) |\sum_{r \mid q} f(q,r) \mu(q/r)|^{2}$, and one can show using the additive large sieve that the assumptions on $\mathcal{A}$ force $\sum_{q=1}^{\infty} g(q)$ to converge, and therefore the first ``big Oh'' term is $o(Qx)$ as $x/Q \rightarrow \infty$. This all extends and refines a (somewhat differently written) asymptotic formula of Hooley~\cite{hooley9}. Vaughan also proved a companion result that is sharper in some cases where $Q\sum_{n \leq x} |a_n|^2 - Qx \sum_{q=1}^{\infty} g(q)$ is very small.

The reader might wonder how ``Criterion U'', or Vaughan's generalised condition, can be sufficient to get an asymptotic for $\sum_{q \leq Q} \sum_{1 \leq a \leq q} \left|\mathcal{A}(x;q,a) - xf(q,(a,q)) \right|^{2}$, since the conditions appear more or less empty when $q \geq \Psi(x)$. However, since we can apply the conditions with $z$ much larger than $x$, they are actually rather strong and, in particular, force $f(q,r)$ to behave well. Thus Hooley~\cite{hooley10} proved an asymptotic assuming only $\mathcal{A}(z;q,a) = z g(q,a) + O_{A}(z/\log^{A}z)$, which at first sight might seem a truly empty condition. In our work, below, we avoid such infinitary conditions and arguments.

We note two particular properties that are forced on $\mathcal{A}$ if Vaughan's conditions are satisfied:
\begin{itemize}
\item $\mathcal{A}$ is dense, and well distributed in short intervals of the form $[z,z+z/\Psi(z)]$, since we have
$$ \mathcal{A}(z) = zf(1,1) + O(z/\Psi(z)) ; $$

\item $\mathcal{A}$ is well distributed in arithmetic progressions to small moduli, e.g. in the sense that
$$ \sum_{r \leq R} \sum_{1 \leq a \leq r} \max_{z \leq x} \left|\mathcal{A}(z;r,a) - \text{Aver}(\mathcal{A},z;r,(a,r)) \right|^{2} \ll R^{2} \max_{z \leq x} \left(\frac{z}{\Psi(z)}\right)^{2}. $$

\end{itemize}
In order for Vaughan's result to provide a genuine variance asymptotic, dominated by $Q\sum_{n \leq x} |a_n|^2 - Qx \sum_{q=1}^{\infty} g(q)$ rather than the ``big Oh'' terms, on a {\em wide range} of $Q$, these must hold for some fairly large $\Psi(z)$.

Vaughan~\cite{vaughangbdh} remarks that one can potentially adapt his arguments, weighting the $a_n$, to treat sparser $\mathcal{A}$, and the reader may consult Vaughan's companion paper~\cite{vaughangbdh1} for a treatment of ``prime like'' sets $\mathcal{A}$ using weights. But it isn't clear that such an adaptation succeeds in general, especially for really sparse $\mathcal{A}$. In the very recent papers \cite{bruv, bruv2}, Br\"udern and Vaughan give special treatments (again using sophisticated circle method type arguments) of cases such as $a_n$ being the representation function of $n$ as a sum of two cubes, or a sum of certain other combinations of powers. They comment~\cite{bruv} that ``A common feature of existing work seems to be that the underlying sequence is fairly dense. ... We are not aware of other examples where a Montgomery--Hooley theorem is known and the growth rate for the expectation [i.e. for $\mathcal{A}(N)$] is as small as $N^{\theta}$ with some $\theta < 1$.'' And in these set-ups one still needs to know that $\mathcal{A}$ is well distributed in short intervals (and in fact in short intervals and progressions simultaneously) in something like the above sense, since this is used to approximate the exponential sum corresponding to $\mathcal{A}$ by suitable known terms. This type of well distribution, with $\Psi(z)$ large, may not be satisfied by many interesting sparse sets. Later we will consider the case where $\mathcal{A}$ is a set of smooth numbers, which is an example of this.

\vspace{12pt}
In this paper we will prove Barban--Davenport--Halberstam type asymptotics that apply to potentially very sparse sets (at the Br\"udern--Vaughan~\cite{bruv,bruv2} level of sparseness), whose behaviour in short intervals is only a tiny bit regular (but which remain well distributed in arithmetic progressions to small moduli). Furthermore, we will obtain quite general results using fairly simple, quick and direct methods, avoiding the use of exponential sums, Tauberian theorems, contour integration, and the other sophisticated machinery often deployed in this area.

\subsection{Statement of results}
We shall make use of the following assumptions, which hopefully seem fairly simple and clean to state:
\begin{itemize}
\item (Progressions Condition) $\mathcal{A} = (a_n)_{n \leq x}$ is well distributed in arithmetic progressions to small moduli, in the sense that
$$ \sum_{1 \leq r \leq 2x/Q} \sum_{1 \leq a \leq r} \Biggl( \sum_{n \leq x, n \equiv a \; \textrm{mod } r} |a_{n}| \Biggr) \max_{z \leq x} \left| \mathcal{A}(z;r,a) - \text{Aver}(\mathcal{A},z;r,(a,r)) \right| \leq \frac{(\sum_{n \leq x} |a_n|)^2}{K_{\text{prog}}} $$
for some $K_{\text{prog}} = K_{\text{prog}}(\mathcal{A},x,Q) > 0$.

\item (Non-concentration on Multiples) for all $H \leq Q$ and for some $K_{\text{conc}}[H] = K_{\text{conc}}(\mathcal{A},x,Q)[H] > 0$, we have
$$ \sum_{H \leq h \leq Q} \Biggl( \sum_{\substack{n \leq x, \\ h \mid n}} |a_{n}| \Biggr)^2 \leq \frac{(\sum_{n \leq x} |a_n|)^2}{K_{\text{conc}}[H]} . $$
\end{itemize}

\vspace{12pt}
Our first result can supply an asymptotic for the variance $V(\mathcal{A},x,Q)$ in many cases where the set of interest, or more generally the support of $(a_n)_{n \leq x}$, is sparse. It may seem strange that we could potentially do better in the sparse case than the dense case, but this is simply because the easy diagonal contribution is quite likely to dominate the variance in the sparse case, where there are fewer other terms around.

\begin{thm1}\label{basicthm}
Suppose that $\sqrt{2x} < Q \leq x$ are large, and that $\mathcal{A}$ satisfies the Progressions Condition and the Non-concentration on Multiples condition. Suppose furthermore that the following holds:
\begin{itemize}
\item (Hereditarily Sparse) for some $K_{\text{hered}} = K_{\text{hered}}(\mathcal{A},x,Q) > 0$, and for any $h \leq 2x/Q$ and any interval $I$ of length $Q/2$, we have
$$ \sum_{n \in I} |a_{nh}| \tau(n) \leq \frac{1}{K_{\text{hered}}} \frac{(\sum_{n \leq x} |a_n|)^2}{\sum_{n \leq x} |a_n| \tau_{3}(n)} , $$
where $\tau(n)$ is the divisor function, and $\tau_{3}(n) := \#\{(u,v,w) \in \N^{3}: uvw=n\}$ is the threefold divisor function.
\end{itemize}

Then $V(\mathcal{A},x,Q)$ is
\begin{eqnarray}
& = & \frac{Q}{2} \sum_{n \leq x} |a_n|^2 + O\left( (\sum_{n \leq x} |a_n|)^2 \left( \frac{1}{K_{\text{prog}}} + \frac{1}{K_{\text{conc}}[2x/Q]} + \frac{\log^{2}(2x/Q)}{K_{\text{hered}}} \right) + \sum_{n \leq x} |a_n|^2 \tau(n) \right) . \nonumber
\end{eqnarray}
\end{thm1}

Note that the Hereditarily Sparse condition demands only upper bound information, rather than asymptotics, for $\mathcal{A}$ in short intervals.

If the sequence $(a_n)_{n \leq x}$ is fairly nicely distributed then we might hope that typically $\sum_{\substack{n \leq x, \\ h \mid n}} |a_{n}| \approx \frac{1}{h} \sum_{n \leq x} |a_n|$, and so we could take $K_{\text{conc}}[2x/Q] \approx x/Q$. Then the error term in Theorem \ref{basicthm} involving $K_{\text{conc}}[2x/Q]$ would be $\approx \frac{Q}{x} (\sum_{n \leq x} |a_n|)^2$. In particular, if $(a_n)_{n \leq x}$ is the indicator function of a set with density $\alpha$ then this error term would be $\approx \alpha^2 Qx$, and so potentially negligible compared with the ``main term'' $\frac{Q}{2} \sum_{n \leq x} |a_n|^2 = \frac{\alpha Qx}{2}$ provided only that we are in the sparse regime where $\alpha$ is small.

Similarly, in the Hereditarily Sparse condition we might hope that up to fairly small factors (e.g. powers of log) arising from the divisor functions $\tau(n), \tau_{3}(n)$, we would generally have $\sum_{n \in I} |a_{nh}| \tau(n) \lesssim \frac{Q}{x} \frac{(\sum_{n \leq x} |a_n|)^2}{\sum_{n \leq x} |a_n| \tau_{3}(n)}$, because we are summing over an interval $I$ of length $Q/2$ rather than the full range up to $x$. Then we could take $K_{\text{hered}} \approx x/Q$ (up to logarithmic type factors, say), and again for $(a_n)_{n \leq x}$ the indicator function of a set with density $\alpha$, the $K_{\text{hered}}$-error term in Theorem \ref{basicthm} might be negligible provided $\alpha$ is at most a certain negative power of $\log x$.

As a very quick illustration, if we take $a_n = \Lambda(n)$ and $x^{0.51} \leq Q \leq x$ (say) then e.g. the classical Bombieri--Vinogradov theorem implies we can take $K_{\text{prog}} \asymp_{A} \log^{A}x$ for any fixed $A > 0$. We may trivially take $K_{\text{conc}}[2x/Q] \asymp x/Q$, and indeed something much stronger. A standard sieve bound also implies that $\sum_{n \in I} |a_{nh}| \tau(n) \ll \log h + \sum_{n \in I \cap [1,x/h]} \Lambda(n) \tau(n) \ll Q$, so we can take $K_{\text{hered}} \asymp x/Q$, and Theorem \ref{basicthm} yields that
$$ V(\Lambda,x,Q) = \frac{Q}{2} \sum_{n \leq x} \Lambda(n)^2 + O_{A}\left( \frac{x^2}{\log^{A}x} + Qx\log^{2}(\frac{2x}{Q}) \right) = \frac{Qx\log x}{2} + O_{A}\left( \frac{x^2}{\log^{A}x} + Qx\log^{2}(\frac{2x}{Q}) \right) . $$
This is of roughly the same quality as the asymptotic first obtained by Montgomery~\cite{mont2} (with a factor $\log^{2}(\frac{2x}{Q})$ rather than $\log(\frac{2x}{Q})$).

\vspace{12pt}
We remark that there is a simple class of sequences $\mathcal{A}$ which shows that imposing conditions on the ``dilates'' $a_{nh}$, as in our Non-concentration on Multiples condition and our Hereditarily Sparse condition, is quite natural when studying $V(\mathcal{A},x,Q)$. Consider the case where $a_n := \textbf{1}_{p|n}$, for some given large prime $p \leq \sqrt{\frac{x}{2}}$, say (and with $\sqrt{2x} < Q \leq x$ as in Theorem \ref{basicthm}). For all moduli $q$ not divisible by $p$, we have $\mathcal{A}(x;q,a) = \N(x/p;q,a\overline{p})$, where $\overline{p}$ is the multiplicative inverse of $p$ modulo $q$. For moduli $q=q'p$ that are divisible by $p$, we have $\mathcal{A}(x;q,a) = 0$ unless $a=a'p$ is also divisible by $p$, in which case $\mathcal{A}(x;q,a) = \N(x/p;q',a')$. This all implies that $V(\mathcal{A},x,Q)$ is
$$ = V(\N,x/p,Q) - \sum_{\substack{Q/2 < q \leq Q, \\ p|q}} \sum_{1 \leq a \leq q} \left|\N(x/p;q,a) - \text{Aver}(\N,x/p;q,(a,q)) \right|^{2} + V(\N,x/p,Q/p) . $$
We can estimate the first and third terms here quite precisely using Lemma \ref{nvarlem}, below, and deduce that $V(\mathcal{A},x,Q) \leq (x/p)^2 \int_{x/Qp}^{2x/Qp} \{v\} (1 - \{v\}) \frac{dv}{v^3} + O((Q/p)^2 + Q\log^{2}Q + x/p)$, where $\{\cdot\}$ denotes the fractional part. In particular, since the {\em larger} integral $(x/p)^2 \int_{x/Qp}^{2x/Qp} v \frac{dv}{v^3} = \frac{Qx}{2p}$ we see that $V(\mathcal{A},x,Q)$ is certainly {\em not} $\sim \frac{Qx}{2p}$ (the main term suggested by Theorem \ref{basicthm}) when $p \ll x/Q$ is large, despite the fact that $K_{\text{prog}}(\mathcal{A},x,Q)$ may be taken to be large here. (A quick calculation, similar to the above, shows we may take $K_{\text{prog}} \asymp Q/p$.) This discrepancy is captured by the terms $K_{\text{conc}}[2x/Q], K_{\text{hered}}$.

We also quickly observe that even if we only know that $|a_n| \leq 1$, we can still use the easy bound $\sum_{2x/Q \leq h \leq Q} ( \sum_{\substack{n \leq x, \\ h \mid n}} |a_{n}| )^2 \leq Q \sum_{2x/Q \leq h \leq Q} \sum_{\substack{n \leq x, \\ h \mid n}} |a_{n}| \leq Q \sum_{n \leq x} |a_{n}| \tau(n)$ in our Non-concentration on Multiples condition, and the easy bound $\sum_{n \in I} |a_{nh}| \tau(n) \leq \sum_{n \in I \cap [1,x/h]} \tau(n) \ll Q\log x$ in our Hereditarily Sparse condition. Estimates of this style also arise, implicitly and explicitly, in the work of Hooley~\cite{hooley10,hooley9} and of Vaughan~\cite{vaughangbdh1,vaughangbdh}. Even using these easy bounds, we obtain a reasonable variance upper bound from Theorem \ref{basicthm}.

\vspace{12pt}
Our second (longer and scarier looking) result can provide a Barban--Davenport--Halberstam type asymptotic for sequences $(a_n)_{n \leq x}$ that are well distributed in progressions to small moduli, and that sufficiently resemble the set of all integers in their behaviour relative to small divisors (i.e. are nearly invariant under ``dilation'' by small integers $h$). This should cover many sets and sequences that are not too sparse, although certainly not all such sets (e.g. those produced by a sieve process, which forbids some small divisors, would not be allowed). It is possible that one could generalise the proof to allow more diverse small divisor behaviour, e.g. resembling other well understood ``model'' sets rather than the set of all integers, but we do not pursue this. 

\begin{thm2}\label{likezthm}
Suppose that $\sqrt{2x} < Q \leq x$ are large, and that $\mathcal{A}$ satisfies the Progressions Condition and Non-concentration on Multiples condition. Suppose furthermore that for some parameters $3 \leq P \leq R \leq x^{1/10}$, and for some $K_{\text{int}}^{(1)}, K_{\text{int}}^{(2)} > 0$ and $K_{\text{int}}^{(3)} \geq 1$ (all possibly depending on $\mathcal{A},x,Q$), the following hold:
\begin{itemize}
\item (Resembles the Integers I) for any $h \in \N$ and any interval $I \subseteq [1,x/h]$, we have
$$ \Biggl| \sum_{\substack{n \in I, \\ p|n \Rightarrow p > P}} a_{nh} - \frac{\sum_{n \leq x} a_n}{x} \sum_{\substack{n \in I, \\ p|n \Rightarrow p > P}} 1 \Biggr| \leq \frac{|I|}{K_{\text{int}}^{(1)}} \frac{\sum_{n \leq x} |a_n|}{x} + \frac{\sum_{n \leq x} |a_n|}{K_{\text{int}}^{(2)}} , $$
where $|I|$ denotes the length of $I$.

\item (Resembles the Integers II) for any $h \in \N$ and interval $I \subseteq [1,x/h]$, we have
$$ \sum_{n \in I} |a_{nh}| \tau(n) \leq K_{\text{int}}^{(3)} |I| \frac{\sum_{n \leq x} |a_n|}{x} + \frac{\sum_{n \leq x} |a_n|}{K_{\text{int}}^{(2)}} . $$
\end{itemize}

Then $V(\mathcal{A},x,Q)$ is
\begin{eqnarray}
& = & \frac{Q}{2} \sum_{n \leq x} |a_n|^2 - \frac{Q}{2} \frac{|\sum_{n \leq x} a_n |^2}{x} + |\sum_{n \leq x} a_n |^2 \int_{x/Q}^{2x/Q} \{v\} (1 - \{v\}) \frac{dv}{v^3} + O\left(\sum_{n \leq x} |a_n|^2 \tau(n) \right) + \nonumber \\
&& + O\left( (\sum_{n \leq x} |a_n|)^2 \left( \frac{1}{K_{\text{prog}}} + \frac{\log x}{P(\log P) K_{\text{conc}}[1]} + \frac{\log P}{K_{\text{int}}^{(1)}} \sqrt{\frac{Q}{x K_{\text{conc}}[2x/Q]}} \right)  \right) + \nonumber \\
&& + O\left( (\sum_{n \leq x} |a_n|) \log^{6}(\frac{Px}{Q}) \left( \frac{Q}{x K_{\text{int}}^{(1)}} \sum_{n \leq x} |a_n| \tau_{3}^{\frac{Px}{Q}}(n) + (\log^{3}x) (\frac{RP}{K_{\text{int}}^{(2)}} + \frac{Q K_{\text{int}}^{(3)}}{x R^{1/\log P}})\sum_{n \leq x} |a_{n}| \tau_{3}(n) \right) \right) , \nonumber
\end{eqnarray}
where $\{\cdot\}$ denotes the fractional part, and $\tau_{3}^{Px/Q}(n)$ is the multiplicative function that agrees with the threefold divisor function $\tau_{3}(n)$ on powers of primes that are $\leq Px/Q$, and is identically 1 on powers of larger primes.
\end{thm2}

Note that the asymptotic potentially provided by Theorem \ref{likezthm} is more precise than in Theorem \ref{basicthm}, with three ``main terms''. If we again think of the case where $(a_n)_{n \leq x}$ is the indicator function of a set with density $\alpha$, these terms become $\frac{\alpha Qx}{2} - \frac{\alpha^2 Qx}{2} + \alpha^2 x^2 \int_{x/Q}^{2x/Q} \{v\} (1 - \{v\}) \frac{dv}{v^3}$. Here $\frac{\alpha Qx}{2} - \frac{\alpha^2 Qx}{2} = \frac{\alpha(1-\alpha)Qx}{2}$ has the familiar shape of a variance for the sum of independent Bernoulli random variables with success probability $\alpha$. The term $\alpha^2 x^2 \int_{x/Q}^{2x/Q} \{v\} (1 - \{v\}) \frac{dv}{v^3}$ is a scaled version of the variance of the set of all natural numbers in arithmetic progressions. As $\alpha$ becomes smaller (i.e. for sparser sets), the terms involving $\alpha^2$ become less significant relative to the simple diagonal contribution $\frac{\alpha Qx}{2}$.

The integral $\int_{x/Q}^{2x/Q} \{v\} (1 - \{v\}) \frac{dv}{v^3}$ could be analysed further (e.g. noting that $\{v\} (1 - \{v\}) = 1/6 - B_{2}(\{v\})$, where $B_2$ is the second Bernoulli polynomial), but seems sufficiently explicit for most purposes. In particular, it is easy to see that this integral is always $\asymp (Q/x)^2$ when $Q \leq x$, and it is $= (1/16)(Q/x)^2 + O((Q/x)^3)$ when $x/Q \rightarrow \infty$.

A key thing to note about Theorem \ref{likezthm} is that it can still provide an asymptotic for $V(\mathcal{A},x,Q)$ even if the relative saving $\frac{1}{K_{\text{int}}^{(1)}}$ assumed in the Resembles the Integers I condition is quite weak. Thus if we have $K_{\text{conc}}[2x/Q] \approx x/Q$ and $\sum_{n \leq x} |a_n| \tau_{3}^{Px/Q}(n) \lesssim \log^{O(1)}(Px/Q) \sum_{n \leq x} |a_n|$, then the error terms involving $K_{\text{int}}^{(1)}$ will both have size $\lesssim (\sum_{n \leq x} |a_n|)^2 \frac{Q}{x} \frac{\log^{O(1)}(Px/Q)}{K_{\text{int}}^{(1)}}$, so it should suffice if $K_{\text{int}}^{(1)}$ is rather larger than $\log^{O(1)}(Px/Q)$ (which might be a power of $\log\log x$, potentially). In particular, this means that Theorem \ref{likezthm} may be applicable even to sequences $\mathcal{A}$ whose counting functions are a bit irregular, and cannot be very well approximated by very simple smooth functions.

\vspace{12pt}
As a new application of both of our Theorems, we deduce a Barban--Davenport--Halberstam type asymptotic for the $y$-smooth numbers less than $x$, on a wide range of $Q$ and $y$. Recall that a number is said to be {\em $y$-smooth} if all of its prime factors are at most $y$. In section 4.2 of his survey~\cite{granvillesmoothsurvey} on smooth numbers, Granville discusses this question, writing ``Bob Vaughan and I have noted that we can get a nontrivial upper bound, but have had difficulties obtaining an asymptotic... for various ranges of values of $y$''. Vaughan has continued to investigate the smooth numbers problem, and we mention his recent work (kindly shared in personal communication with the author) below.

To describe our result, we introduce the standard notation for counting functions of the $y$-smooth numbers: letting $\mathcal{S}(y)$ be the set of all $y$-smooth numbers, we put
$$ \Psi(x,y) := \sum_{n \leq x} \textbf{1}_{\{n \in \mathcal{S}(y)\}}, \;\;\; \Psi(x,y;q,a) := \sum_{\substack{n \leq x, \\ n \equiv a (\textrm{mod } q)}} \textbf{1}_{\{n \in \mathcal{S}(y)\}}, \;\;\; \Psi_q(x,y) := \sum_{\substack{n \leq x, \\ (n,q)=1}} \textbf{1}_{\{n \in \mathcal{S}(y)\}} . $$
Note that in this case we have $\Psi(x,y;q,a) = \Psi(x/h,y;q/h,a/h)$ if $h = (a,q)$ is $y$-smooth, and we have $\Psi(x,y;q,a) = 0$ otherwise, so the Barban--Davenport--Halberstam type variance of the smooth numbers can also be written as
$$ V(\mathcal{S}(y),x,Q) = \sum_{Q/2 < q \leq Q} \sum_{\substack{h \mid q, \\ h \in \mathcal{S}(y)}} \sum_{\substack{1 \leq a \leq q, \\ (a,q) = h}} \left|\Psi(x,y;q,a) - \frac{\Psi_{q/h}(x/h,y)}{\phi(q/h)} \right|^{2} . $$

\begin{cor1}\label{smoothcor}
There exist a small absolute constant $c > 0$, and a large absolute constant $C > 0$, such that the following is true.

Suppose that $\log^{C}x \leq y \leq x$ are large, and $x^{0.51} \leq Q \leq x$. For any $A > 0$, we have
\begin{eqnarray}
V(\mathcal{S}(y),x,Q) & = & \frac{Q}{2} \Psi(x,y) \left(1 - \frac{\Psi(x,y)}{x}\right) + \Psi(x,y)^2 \int_{x/Q}^{2x/Q} \{v\} (1 - \{v\}) \frac{dv}{v^3} + \nonumber \\
&& + O_{A}\left( \Psi(x,y)^2 \left( \frac{e^{-\frac{cu}{\log^{2}(u+1)}}}{\log^{A}x} + y^{-c} \right) + Q\Psi(x,y) \left( \frac{e^{-cu\log(u+1)} (\log\log x)^{12}}{\log x} \right) \right) . \nonumber
\end{eqnarray}
Here $u := (\log x)/\log y$, and the ``big Oh'' implicit constant may depend on $A$.
\end{cor1}

One has $\Psi(x,y) \asymp x$ whenever $x^{0.01} \leq y \leq x$, say, whereas $\Psi(x,\log^{C}x) = x^{1 - 1/C + o(1)}$ (see e.g. Corollary 7.9 of Montgomery and Vaughan~\cite{mv}). Thus Corollary \ref{smoothcor} can handle levels of sparseness comparable to those discussed by Br\"udern and Vaughan~\cite{bruv}.

Corollary \ref{smoothcor} provides a true asymptotic for $V(\mathcal{S}(y),x,Q)$, with the ``main terms'' on the first line dominating the ``big Oh'' terms, when $\log^{C}x \leq y \leq \frac{x}{e^{(\log\log x)^{13}}}$ (say, so that $\left(1 - \frac{\Psi(x,y)}{x}\right) \gg \frac{(\log\log x)^{13}}{\log x}$) and $Q$ is larger than $\Psi(x,y) (\frac{e^{-\frac{cu}{\log^{2}(u+1)}}}{\log^{A}x} + y^{-c})$; and also when $\frac{x}{e^{(\log\log x)^{13}}} < y \leq x$, on the smaller range $x \frac{(\log\log x)^{12}}{\log x} \ll Q \leq x$ where the integral involving $v$ dominates. When $y$ is so large, $\mathcal{S}(y)$ is very close to simply being the set of all natural numbers (only the small subset having one very large prime factor $> y$ are excluded), so if desired one could probably refine Corollary \ref{smoothcor} to extend this range of $Q$ by working more carefully.

Vaughan (personal communication) has recently obtained an asymptotic for a quantity like $V(\mathcal{S}(y),x,Q)$, with the subtracted terms $\frac{\Psi_{q/h}(x/h,y)}{\phi(q/h)}$ replaced by a more complicated expression suggested by the circle method. His unconditional result is valid for $e^{(\log\log x)^{5/3 + \epsilon}} \leq y \leq x^{\frac{2}{3(\log\log x)^4}}$ and $Q \leq \Psi(x,y)$, and has error terms generally a little worse than in Corollary \ref{smoothcor}. He also formulates a stronger result (valid for $e^{(\log\log x)^{5/3 + \epsilon}} \leq y \leq x^{1/\log\log x}$, and with some improved error terms), under a suitable hypothesis that Dirichlet $L$-functions have no exceptional (Siegel) zeros.

This may all be compared with the known upper and lower bounds for $V(\mathcal{S}(y),x,Q)$. Harper~\cite{harpersmoothbv} proved an upper bound $\ll_{A} Q\Psi(x,y) + \Psi(x,y)^2 (\frac{e^{-\frac{cu}{\log^{2}(u+1)}}}{\log^{A}x} + y^{-c})$ for (a slightly different version of) the variance, on the range $\log^{C}x \leq y \leq x$ and $Q \leq \Psi(x,y)$. When $\sqrt{x} \leq y \leq x/C$ and $x^{1/2 + \delta} \leq Q \leq x$, Mastrostefano~\cite{mast} proved the lower bound $\gg_{\delta} Qx\log u + Q^2$, again for a slightly different version of the variance. 

% SECTION 2 %%%%%%%%%%%%%%%%%%%%%%%%%%%%%
\section{Preliminary manoeuvres}\label{prelimsec}
We may assume throughout that $\mathcal{A} = (a_n)_{n \leq x}$ is a {\em real} sequence, since in the general complex case we have $V(\mathcal{A},x,Q) = V(\Re\mathcal{A},x,Q) + V(\Im\mathcal{A},x,Q)$ (with the obvious notation $\Re\mathcal{A} = (\Re(a_n))_{n \leq x}$ and $\Im\mathcal{A} = (\Im(a_n))_{n \leq x}$), and the reader may check that the main and error terms in Theorems \ref{basicthm} and \ref{likezthm} then sum in the desired way.

In this section we will derive a preliminary expression for our variance $V(\mathcal{A},x,Q)$, that will be used in the proofs of both of our theorems. In doing this we will make use (in a standard way) of our Progressions Condition, which will not then be used again. After that we will reorganise the terms in an elementary but non-standard way, which will be important in the following sections to avoid incurring unacceptable error terms.

\vspace{12pt}
We begin with some straightforward manipulations, noting that $V(\mathcal{A},x,Q)$ is
\begin{eqnarray}\label{stopensquare}
& = &  \sum_{Q/2 < q \leq Q} \sum_{h \mid q} \sum_{1 \leq a \leq q, (a,q)=h} \left|\mathcal{A}(x;q,a) - \text{Aver}(\mathcal{A},x;q,h) \right|^{2} \nonumber \\
& = & \sum_{Q/2 < q \leq Q} \sum_{1 \leq a \leq q} \mathcal{A}(x;q,a)^{2} - \sum_{Q/2 < q \leq Q} \sum_{h \mid q} \phi(q/h) \text{Aver}(\mathcal{A},x;q,h)^{2} \nonumber \\
& = & \sum_{Q/2 < q \leq Q} \sum_{1 \leq a \leq q} \Biggl( \sum_{\substack{n \leq x, \\ n \equiv a \; \textrm{mod } q}} a_n \Biggr)^{2} - \sum_{Q/2 < q \leq Q} \sum_{h \mid q} \frac{1}{\phi(q/h)} \Biggl( \sum_{n \leq x, (n,q)=h} a_n \Biggr)^{2} .
\end{eqnarray}
We can rewrite the first sums in \eqref{stopensquare} as
\begin{eqnarray}
&& \sum_{Q/2 < q \leq Q} \sum_{n_1 , n_2 \leq x} a_{n_1} a_{n_2} \textbf{1}_{q \mid n_1 - n_2} = (\frac{Q}{2} + O(1))\sum_{n \leq x} a_n^2 +  \sum_{\substack{n_1 , n_2 \leq x, \\ n_1 \neq n_2}} a_{n_1} a_{n_2} \sum_{Q/2 < q \leq Q} \textbf{1}_{n_1 - n_2 = qr \; \textrm{for some } r} \nonumber \\
& = & (Q/2 + O(1)) \sum_{n \leq x} a_n^2 + 2 \sum_{n_2 < n_1 \leq x} a_{n_1} a_{n_2} \sum_{Q/2 < q \leq Q} \textbf{1}_{n_1 - n_2 = qr \; \textrm{for some } r} , \nonumber
\end{eqnarray}
and then apply the ingenious {\em divisor switching idea} of Hooley~\cite{hooley}, whereby we switch from a summation over $q$ to a summation over the complementary divisor $r$ (which must now be a positive integer of size at most $2x/Q$). Thus the previous line is
$$ = (Q/2 + O(1))\sum_{n \leq x} a_n^2 + 2 \sum_{1 \leq r \leq 2x/Q} \sum_{1 \leq a \leq r} \sum_{\substack{n_1 \leq x, \\ n_1 \equiv a \; \textrm{mod } r}} a_{n_1} \sum_{\substack{n_2 \leq x, \\ n_2 \equiv a \; \textrm{mod } r}} a_{n_2} \textbf{1}_{Q/2 < (n_1 - n_2)/r \leq Q} . $$
The point of this manoeuvre is that all of our congruence conditions are now to small moduli, so we can use our Progressions Condition. Indeed, the previous line is
\begin{eqnarray}
& = & (\frac{Q}{2} + O(1))\sum_{n \leq x} a_n^2 + 2 \sum_{1 \leq r \leq 2x/Q} \sum_{h \mid r} \sum_{\substack{1 \leq a \leq r, \\ (a,r)=h}} \sum_{\substack{n_1 \leq x, \\ n_1 \equiv a \; \textrm{mod } r}} a_{n_1} \frac{1}{\phi(r/h)} \sum_{\substack{n_2 \leq x, \\ (n_2,r) = h}} a_{n_2} \textbf{1}_{Q/2 < (n_1 - n_2)/r \leq Q} + \nonumber \\
&& + O\Biggl( \sum_{1 \leq r \leq 2x/Q} \sum_{h \mid r} \sum_{\substack{1 \leq a \leq r, \\ (a,r)=h}} \sum_{\substack{n_1 \leq x, \\ n_1 \equiv a \; \textrm{mod } r}} |a_{n_1}| \max_{z \leq x} \left| \mathcal{A}(z;r,a) - \text{Aver}(\mathcal{A},z;r,h) \right| \Biggr) , \nonumber
\end{eqnarray}
and the ``big Oh'' terms here are $\ll \sum_{n \leq x} a_n^2 + \frac{(\sum_{n \leq x} |a_n|)^2}{K_{\text{prog}}}$, which is acceptable for both Theorems \ref{basicthm} and \ref{likezthm}. We can also slightly rewrite the multiple sums, as
\begin{equation}\label{afterprogcond}
2 \sum_{1 \leq r \leq 2x/Q} \sum_{h \mid r} \frac{1}{\phi(r/h)} \sum_{\substack{n_1 \leq x, \\ (n_1,r)=h}} a_{n_1} \sum_{\substack{n_2 \leq x, \\ (n_2,r) = h}} a_{n_2} \textbf{1}_{Q/2 < (n_1 - n_2)/r \leq Q} .
\end{equation}

\vspace{12pt}
At this point, the standard approach would probably be to try expressing the inner sums in \eqref{afterprogcond} in terms of ``nice'' smooth functions. Instead, our philosophy will be to keep all such sums unevaluated until the end, reorganising them to find cancellation.

To explain further, let us consider the second term in our expansion \eqref{stopensquare} of $V(\mathcal{A},x,Q)$, which is untouched so far. That term $\sum_{Q/2 < q \leq Q} \sum_{h \mid q} \frac{1}{\phi(q/h)} ( \sum_{n \leq x, (n,q)=h} a_n )^{2}$ is
$$ = 2  \sum_{Q/2 < q \leq Q} \sum_{h \mid q} \frac{1}{\phi(q/h)} \sum_{\substack{n_2 < n_1 \leq x, \\ (n_1,q)=(n_2,q)=h}} a_{n_1} a_{n_2} + \sum_{Q/2 < q \leq Q} \sum_{h \mid q} \frac{1}{\phi(q/h)} \sum_{n \leq x, (n,q)=h} a_n^2 , $$
and the second sum here is at most
$$ \sum_{n \leq x} a_n^2 \sum_{Q/2 < q \leq Q} \sum_{h \mid (n,q)} \frac{1}{\phi(q/h)} = \sum_{n \leq x} a_n^2 \sum_{h \mid n} \sum_{Q/2h < q' \leq Q/h} \frac{1}{\phi(q')} . $$
It is easy to show the uniform bound $\sum_{Q/2h < q' \leq Q/h} \frac{1}{\phi(q')} \ll 1$ (see e.g. chapter 2.1 of Montgomery and Vaughan~\cite{mv} for much more precise statements), so we get an acceptable overall contribution of $O(\sum_{n \leq x} a_n^2 \tau(n))$ from this second sum, where $\tau(n)$ is the divisor function. And the form of the term $2  \sum_{Q/2 < q \leq Q} \sum_{h \mid q} \frac{1}{\phi(q/h)} \sum_{\substack{n_2 < n_1 \leq x, \\ (n_1,q)=(n_2,q)=h}} a_{n_1} a_{n_2}$ is very like the form of the term we were working with previously in \eqref{afterprogcond}, where $q$ is replaced by $r$ (and some of the ranges of summation are slightly different).

We can make this even clearer if we reorganise a little more, noting that \eqref{afterprogcond} is
\begin{eqnarray}
& = & 2 \sum_{h \leq 2x/Q} \sum_{d \leq 2x/Qh} \frac{1}{\phi(d)} \sum_{\substack{n_1 \leq x, \\ (n_1,dh)=h}} a_{n_1} \sum_{\substack{n_2 \leq x, \\ (n_2,dh) = h}} a_{n_2} \textbf{1}_{Q/2 < \frac{n_1 - n_2}{dh} \leq Q} \nonumber \\
& = & 2 \sum_{h \leq 2x/Q} \sum_{\substack{n_2 < n_1 \leq x, \\ h \mid n_1, h \mid n_2}} a_{n_1} a_{n_2} \sum_{\substack{(n_1 - n_2)/(Qh) \leq d < 2(n_1 - n_2)/(Qh), \\ (n_1 n_2 /h^2,d)=1}} \frac{1}{\phi(d)} , \nonumber
\end{eqnarray}
and similarly
$$ 2  \sum_{Q/2 < q \leq Q} \sum_{h \mid q} \frac{1}{\phi(q/h)} \sum_{\substack{n_2 < n_1 \leq x, \\ (n_1,q)=(n_2,q)=h}} a_{n_1} a_{n_2} = 2 \sum_{h \leq Q} \sum_{\substack{n_2 < n_1 \leq x, \\ h \mid n_1, h \mid n_2}} a_{n_1} a_{n_2} \sum_{\substack{Q/(2h) < d \leq Q/h, \\ (n_1 n_2 /h^2,d)=1}} \frac{1}{\phi(d)} . $$

\vspace{12pt}
Putting everything together, and noting that $2x/Q < Q$ if we assume that $Q > \sqrt{2x}$, we have shown in this section that $V(\mathcal{A},x,Q)$ is
\begin{eqnarray}
& = & (Q/2)\sum_{n \leq x} a_n^2 + 2 \sum_{h \leq 2x/Q} \sum_{\substack{n_2 < n_1 \leq x, \\ h \mid n_1, h \mid n_2}} a_{n_1} a_{n_2} \Biggl( \sum_{\substack{\frac{(n_1 - n_2)}{Qh} \leq d < \frac{2(n_1 - n_2)}{Qh}, \\ (n_1 n_2 /h^2,d)=1}} \frac{1}{\phi(d)} - \sum_{\substack{Q/(2h) < d \leq Q/h, \\ (n_1 n_2 /h^2,d)=1}} \frac{1}{\phi(d)} \Biggr) - \nonumber \\
&& - 2 \sum_{2x/Q < h \leq Q} \sum_{\substack{n_2 < n_1 \leq x, \\ h \mid n_1, h \mid n_2}} a_{n_1} a_{n_2} \sum_{\substack{Q/(2h) < d \leq Q/h, \\ (n_1 n_2 /h^2,d)=1}} \frac{1}{\phi(d)} + O\Biggl( \sum_{n \leq x} a_n^2 \tau(n) + \frac{(\sum_{n \leq x} |a_n|)^2}{K_{\text{prog}}} \Biggr) . \nonumber
\end{eqnarray}

For given $x,Q$, let us now define
$$ \mathcal{S}_{h,l} := \{(n_1,n_2) : 1 \leq n_2 < n_1 \leq x/h, \;\; \text{and} \;\; \frac{lQ}{2} < n_1 - n_2 \leq \frac{(l+1)Q}{2} \}, \;\;\;\;\; h \in \N, \; l \in \N \cup \{0\}, $$
and also
$$ C_{m,l,z} := \sum_{\substack{l/2 \leq d < l, \\ (d,m)=1}} \frac{1}{\phi(d)} - \sum_{\substack{z/2 < d \leq z, \\ (d,m)=1}} \frac{1}{\phi(d)}, \;\;\;\;\; m \in \N, \;\;\; l,z > 0. $$
Then we can express our conclusion more compactly: we have
\begin{eqnarray}\label{keyexpression}
V(\mathcal{A},x,Q) & = & \frac{Q}{2} \sum_{n \leq x} a_n^2 + 2 \sum_{h \leq 2x/Q} \sum_{0 \leq l \leq 2x/Qh} \sum_{(n_1,n_2) \in \mathcal{S}_{h,l}} a_{n_1 h} a_{n_2 h} C_{n_1 n_2 , l+1, Q/h} - \nonumber \\
&& - 2 \sum_{2x/Q < h \leq Q} \sum_{\substack{n_2 < n_1 \leq x, \\ h \mid n_1, h \mid n_2}} a_{n_1} a_{n_2} \sum_{\substack{Q/(2h) < d \leq Q/h, \\ (n_1 n_2 /h^2,d)=1}} \frac{1}{\phi(d)} + \text{error},
\end{eqnarray}
where the ``error'' (bounded as described above) is acceptably small for our Theorems.

% SECTION 3 %%%%%%%%%%%%%%%%%%%%%%%%%%%%%
\section{Theorem \ref{basicthm}: A simple result for sparse sets}
To prove Theorem \ref{basicthm}, let us immediately take another look at \eqref{keyexpression}. Again using the easy bound $\sum_{\substack{Q/(2h) < d \leq Q/h, \\ (n_1 n_2 /h^2,d)=1}} \frac{1}{\phi(d)} \leq \sum_{Q/(2h) < d \leq Q/h} \frac{1}{\phi(d)} \ll 1$, and applying the Non-concentration on Multiples condition, we find the sums on the second line of \eqref{keyexpression} are
$$ \ll \sum_{2x/Q < h \leq Q} \Biggl( \sum_{\substack{n \leq x, \\ h \mid n}} |a_{n}| \Biggr)^2 \ll \frac{(\sum_{n \leq x} |a_n|)^2}{K_{\text{conc}}[2x/Q]} . $$

To handle the first line of \eqref{keyexpression} we shall apply the following lemma, which gives some more precise information about the terms $C_{m,l,z}$ occurring there.

\begin{lem1}\label{totientlem}
Suppose that $m \in \N$ and $z \geq l > 0$, and let $C_{m,l,z}$ be as defined above. Then
$$ |C_{m,l,z}| \ll \min\{1,\tau(m) \frac{\log(l+2)}{(l+2)}\}, $$
where $\tau(m)$ denotes the divisor function.
\end{lem1}

\begin{proof}[Proof of Lemma \ref{totientlem}]
This follows easily from the fact (see e.g. Exercise 2.1.1.16 of Montgomery and Vaughan~\cite{mv}) that $\sum_{n \leq z, (n,m)=1} 1/\phi(n) = A_{m}\log z + B_{m} + O(\tau(m)(\log z)/z)$ if $z \geq 2$, for certain constants $A_{m}, B_{m}$ that depend on $m$ only.
\end{proof}

Thus the summands on the first line of \eqref{keyexpression} are typically rather small (because of the way we organised things to form the $C_{m,l,z}$), so we have a good chance of controlling the sums successfully. Indeed, using the sub-multiplicativity of $\tau$ we find that
$$ \Biggl| \sum_{h \leq \frac{2x}{Q}} \sum_{0 \leq l \leq \frac{2x}{Qh}} \sum_{(n_1,n_2) \in \mathcal{S}_{h,l}} a_{n_1 h} a_{n_2 h} C_{n_1 n_2 , l+1, Q/h} \Biggr| \ll \sum_{h \leq \frac{2x}{Q}} \sum_{0 \leq l \leq \frac{2x}{Qh}} \frac{\log(l+2)}{l+2} \sum_{(n_1,n_2) \in \mathcal{S}_{h,l}} |a_{n_1 h}| |a_{n_2 h}| \tau(n_1) \tau(n_2) . $$
Recalling the definition of $\mathcal{S}_{h,l}$, and applying our Hereditarily Sparse condition to the sum over the $n_2$ variable, we find this is all
$$ \ll \frac{1}{K_{\text{hered}}} \frac{(\sum_{n \leq x} |a_n|)^2}{\sum_{n \leq x} |a_n| \tau_{3}(n)} \cdot \sum_{h \leq \frac{2x}{Q}} \sum_{0 \leq l \leq \frac{2x}{Qh}} \frac{\log(l+2)}{l+2} \sum_{n_1 \leq x/h} |a_{n_1 h}| \tau(n_1) . $$
Letting $n = n_1 h$, and moving the sum over $h$ to the inside, we can bound everything by
$$ \frac{1}{K_{\text{hered}}} \frac{(\sum_{n \leq x} |a_n|)^2}{\sum_{n \leq x} |a_n| \tau_{3}(n)} \cdot \sum_{0 \leq l \leq \frac{2x}{Q}} \frac{\log(l+2)}{l+2} \sum_{n \leq x} |a_{n}| \tau_{3}(n) \ll \frac{\log^{2}(2 + x/Q) (\sum_{n \leq x} |a_n|)^2}{K_{\text{hered}}} , $$
which finishes the proof of Theorem \ref{basicthm}.
\qed

% SECTION 4 %%%%%%%%%%%%%%%%%%%%%%%%%%%%%
\section{Theorem \ref{likezthm}: Comparing with the integers}
To prove Theorem \ref{likezthm}, the key observation is that $C_{n_1 n_2 , l+1, Q/h}, \sum_{\substack{Q/(2h) < d \leq Q/h, \\ (n_1 n_2 /h^2,d)=1}} \frac{1}{\phi(d)}$ on the right hand side of \eqref{keyexpression} essentially only depend on $n_1, n_2$ via their small prime factors. This is made precise in the following result.

\begin{lem2}\label{totientprimeslem}
For any $m, n \in \N$ and any $l > 0$, we have
$$ |\sum_{\substack{l/2 < d \leq l, \\ (d,mn)=1}} \frac{1}{\phi(d)} -  \sum_{\substack{l/2 < d \leq l, \\ (d,m)=1}} \frac{1}{\phi(d)}| \ll \sum_{p \mid n} \frac{1}{p} , $$
where $p$ denotes primes.

In particular, we always have $|C_{mn, l,z} - C_{m,l,z}| \ll \sum_{p \mid n} \frac{1}{p}$.
\end{lem2}

\begin{proof}[Proof of Lemma \ref{totientprimeslem}]
The left hand side in the first statement is $\leq \sum_{p \mid n} \sum_{\substack{l/2 < d \leq l, \\ p|d}} \frac{1}{\phi(d)} \ll \sum_{p \mid n} \frac{1}{p} \sum_{l/2p < d \leq l/p} \frac{1}{\phi(d)}$, and then we just apply the easy bound $\ll 1$ for the sums over $d$. The second statement follows immediately from the first and the definition of $C_{m,l,z}$.
\end{proof}

So let us temporarily write $\text{Sm}_{P}(n)$ to denote the largest $P$-smooth divisor of $n$, and (recalling our notation from the end of section \ref{prelimsec}) let $\mathcal{S}_{h,l}^{(m_1, m_2)}$ denote the subset of $\mathcal{S}_{h,l}$ consisting of pairs $(n_1,n_2)$ with $\text{Sm}_{P}(n_1) = m_1$ and $\text{Sm}_{P}(n_2) = m_2$. Then Lemma \ref{totientprimeslem} implies that $C_{n_1 n_2 , l+1, Q/h} = C_{m_1 m_2 , l+1, Q/h} + O(\sum_{\substack{p| n_1 n_2, \\ p > P}} \frac{1}{p})$ when $(n_1,n_2) \in \mathcal{S}_{h,l}^{(m_1, m_2)}$, similarly for $\sum_{\substack{Q/(2h) < d \leq Q/h, \\ (n_1 n_2 /h^2,d)=1}} \frac{1}{\phi(d)}$, and we can rewrite the double sums of \eqref{keyexpression} as
\begin{eqnarray}\label{keyexprreplace}
&& 2 \sum_{h \leq 2x/Q} \sum_{0 \leq l \leq 2x/Qh} \sum_{\substack{m_1, m_2, \\ P \; \text{smooth}}} C_{m_1 m_2 , l+1, Q/h} \sum_{(n_1,n_2) \in \mathcal{S}_{h,l}^{(m_1, m_2)}} a_{n_1 h} a_{n_2 h} - \\
&& - 2 \sum_{\frac{2x}{Q} < h \leq Q} \sum_{\substack{m_1, m_2, \\ P \; \text{smooth}}} \sum_{\substack{\frac{Q}{2h} < d \leq \frac{Q}{h}, \\ (m_1 m_2, d)=1}} \frac{1}{\phi(d)} \sum_{\substack{n_2 < n_1 \leq x/h, \\ \text{Sm}_{P}(n_1) = m_1, \\ \text{Sm}_{P}(n_2) = m_2}} a_{n_1 h} a_{n_2 h} + O\Biggl( \sum_{h \leq Q} \sum_{n_2 < n_1 \leq \frac{x}{h}} |a_{n_1 h}| |a_{n_2 h}| \sum_{\substack{p| n_1 n_2, \\ p > P}} \frac{1}{p} \Biggr) . \nonumber
\end{eqnarray}
Proceeding a bit crudely, since $n_1 n_2$ can have at most $\frac{\log(n_1 n_2)}{\log P} \ll \frac{\log x}{\log P}$ prime factors larger than $P$, and using the Non-concentration on Multiples condition, this ``big Oh'' term is $\ll \frac{\log x}{P \log P} \sum_{h \leq Q} \sum_{n_2 < n_1 \leq x/h} |a_{n_1 h}| |a_{n_2 h}| \ll \frac{\log x}{P \log P} \frac{(\sum_{n \leq x} |a_n|)^2}{K_{\text{conc}}[1]}$.

The sums over $m_1, m_2$ in \eqref{keyexprreplace} currently contain too many terms to be easily manageable, so we wish to truncate them to the ranges $m_1, m_2 \leq R$. To do this, we may note that if $\text{Sm}_{P}(n) > R$ then $n$ must have at least one $P$-smooth divisor from the interval $(R,RP]$ (e.g. one can successively remove the largest prime factors from $\text{Sm}_{P}(n)$ until no more can be removed without its size dropping below $R$). So using Lemma \ref{totientlem} along with the easy estimate $\sum_{\substack{\frac{Q}{2h} < d \leq \frac{Q}{h}, \\ (m_1 m_2, d)=1}} \frac{1}{\phi(d)} \ll 1$, we can bound the total contribution from any $m_1 > R$ or $m_2 > R$ by
$$ \sum_{\substack{R < m \leq RP, \\ P \; \text{smooth}}} \Biggl( \sum_{h \leq \frac{2x}{Q}} \sum_{0 \leq l \leq \frac{2x}{Qh}} \frac{\log(l+2)}{l+2} \sum_{\substack{(n_1,n_2) \in \mathcal{S}_{h,l}, \\ m|n_1 \; \text{or} \; m|n_2}} \tau(n_1 n_2) |a_{n_1 h}| |a_{n_2 h}| + \sum_{\frac{2x}{Q} < h \leq Q} \sum_{n_1 \leq x/h} |a_{n_1 h}| \sum_{\substack{n_2 \leq x/h, \\ m|n_2}} |a_{n_2 h}| \Biggr) . $$
Rewriting sums like $\sum_{\substack{n_2 \leq x/h, \\ m|n_2}} |a_{n_2 h}|$ as $\sum_{n \leq x/hm} |a_{n hm}|$, and using the sub-multiplicativity of $\tau$, the Resembles the Integers II condition implies this all has order at most
$$ \sum_{\substack{R < m \leq RP, \\ P \; \text{smooth}}} \Biggl( \sum_{h \leq \frac{2x}{Q}} \log^{2}(2+\frac{x}{Q}) \sum_{n \leq \frac{x}{h}} |a_{n h}| \tau(n m)  (\frac{K_{\text{int}}^{(3)} Q}{mx} + \frac{1}{K_{\text{int}}^{(2)}}) + \sum_{\frac{2x}{Q} < h \leq Q} \sum_{n_1 \leq \frac{x}{h}} |a_{n_1 h}| (\frac{K_{\text{int}}^{(3)}}{hm} + \frac{1}{K_{\text{int}}^{(2)}}) \Biggr) (\sum_{n \leq x} |a_n|) . $$

The contribution from the $\frac{1}{K_{\text{int}}^{(2)}}$ terms here is, again somewhat crudely,
$$ \ll \frac{(\sum_{n \leq x} |a_n|) \log^{2}(2+\frac{x}{Q})}{K_{\text{int}}^{(2)}} (\sum_{\substack{R < m \leq RP, \\ P \; \text{smooth}}} \tau(m)) \Biggl( \sum_{h \leq \frac{2x}{Q}} \sum_{n \leq \frac{x}{h}} |a_{n h}| \tau(n) + \sum_{\frac{2x}{Q} < h \leq Q} \sum_{n_1 \leq \frac{x}{h}} |a_{n_1 h}| \Biggr) . $$
Since $\sum_{\substack{m \leq RP , \\ P \; \text{smooth}}} \tau(m) \ll RP \log P$, we find this is all $\ll \frac{(\sum_{n \leq x} |a_n|) \log^{3}(P+\frac{x}{Q}) RP}{K_{\text{int}}^{(2)}} \sum_{n \leq x} |a_n| \tau_{3}(n)$. Meanwhile, the bound $\sum_{\substack{R < m \leq RP , \\ P \; \text{smooth}}} \frac{\tau(m)}{m} \leq R^{-1/\log P} \sum_{\substack{m : \\ P \; \text{smooth}}} \frac{\tau(m)}{m^{1-1/\log P}} \ll R^{-1/\log P} \log^{2}P$, together with the fact that $1/h \ll Q/x$ when $h > 2x/Q$, imply that the contribution from all the $K_{\text{int}}^{(3)}$ terms is
$$ \ll \frac{K_{\text{int}}^{(3)} Q}{x} (\sum_{n \leq x} |a_n|) \frac{\log^{2}P}{R^{1/\log P}} \Biggl( \log^{2}(2+\frac{x}{Q}) \sum_{h \leq \frac{2x}{Q}} \sum_{n \leq \frac{x}{h}} |a_{n h}| \tau(n) + \sum_{\frac{2x}{Q} < h \leq Q} \sum_{n_1 \leq \frac{x}{h}} |a_{n_1 h}| \Biggr) , $$
which is certainly $\ll \frac{K_{\text{int}}^{(3)} Q}{x} (\sum_{n \leq x} |a_n|) \frac{\log^{4}(P+\frac{x}{Q})}{R^{1/\log P}} \sum_{n \leq x} |a_n| \tau_{3}(n)$. All of these error estimates are acceptable for Theorem \ref{likezthm}.

\vspace{12pt}
Now it remains to analyse the contribution to \eqref{keyexprreplace} when $m_1, m_2 \leq R$. For given $n_1$, the sum of $a_{n_2 h}$ over all $n_2$ such that $(n_1,n_2) \in \mathcal{S}_{h,l}^{(m_1, m_2)}$ may be written as a sum of $a_{n m_2 h}$ over all $n = n_{2}/m_{2}$ belonging to a certain sub-interval of $[1,x/hm_2]$ (of length at most $Q/2m_2$), and which have no prime factors $\leq P$. Similarly, the sum of $a_{n_2 h}$ over all $n_2 < n_1$ with $\text{Sm}_{P}(n_2) = m_2$ is a sum of $a_{n m_2 h}$ over an interval of length at most $x/hm_2$, with a restriction of no prime factors $\leq P$. So using our Resembles the Integers I condition (with $h$ replaced by $hm_2$), we may replace $a_{n_2 h}$ everywhere by the constant $\frac{\sum_{n \leq x} a_n}{x}$, at the cost of an error term that is
\begin{eqnarray}
& \ll & (\sum_{n \leq x} |a_n|) \Biggl( \sum_{h \leq \frac{2x}{Q}} \sum_{0 \leq l \leq \frac{2x}{Qh}} \frac{\log(l+2)}{l+2} \sum_{\substack{m_2 \leq R, \\ P \; \text{smooth}}} \tau(m_2) \sum_{n_1 \leq x/h} \tau(\text{Sm}_{P}(n_1)) |a_{n_1 h}| (\frac{Q}{K_{\text{int}}^{(1)} x m_2} + \frac{1}{K_{\text{int}}^{(2)}} ) + \nonumber \\
&& + \sum_{2x/Q < h \leq Q} \sum_{\substack{m_2 \leq R , \\ P \; \text{smooth}}} \sum_{n_1 \leq x/h} |a_{n_1 h}| (\frac{1}{K_{\text{int}}^{(1)} hm_2} + \frac{1}{K_{\text{int}}^{(2)}} ) \Biggr) . \nonumber
\end{eqnarray}

Arguing exactly as we did above, using the bounds $\sum_{\substack{m_2 \leq R , \\ P \; \text{smooth}}} 1 \leq \sum_{\substack{m_2 \leq R , \\ P \; \text{smooth}}} \tau(m_2) \ll R \log P$ and the fact that $\tau(\text{Sm}_{P}(n_1)) \leq \tau(n_1)$, we find the contribution from all these $\frac{1}{K_{\text{int}}^{(2)}}$ terms is $\ll \frac{(\sum_{n \leq x} |a_n|) R \log^{3}(P + \frac{x}{Q})}{K_{\text{int}}^{(2)}} \sum_{n \leq x} |a_{n}| \tau_{3}(n)$.

Since $\sum_{\substack{m_2 \leq R , \\ P \; \text{smooth}}} \frac{1}{m_2} \leq \sum_{\substack{m_2 : \\ P \; \text{smooth}}} \frac{1}{m_2} \ll \log P$ and $\sum_{\substack{m_2 \leq R , \\ P \; \text{smooth}}} \frac{\tau(m_2)}{m_2} \ll \log^{2}P$, the contribution from all the $\frac{1}{K_{\text{int}}^{(1)}}$ terms is
$$ \ll \frac{(\sum_{n \leq x} |a_n|) }{K_{\text{int}}^{(1)}} \Biggl( \frac{Q}{x} \log^{2}(2+\frac{x}{Q}) \log^{2}P \sum_{h \leq \frac{2x}{Q}} \sum_{n_1 \leq x/h} \tau(\text{Sm}_{P}(n_1)) |a_{n_1 h}| + \log P \sum_{2x/Q < h \leq Q} \frac{1}{h} \sum_{n_1 \leq x/h} |a_{n_1 h}| \Biggr) . $$
When bounding this we shall argue a little more delicately, since we envisage that $K_{\text{int}}^{(1)}$ may not be particularly large in practice. The first sums inside the bracket are $\sum_{n \leq x} |a_n| \sum_{h|n, h \leq 2x/Q} \tau(\text{Sm}_{P}(n/h))$, and since a divisor $h \leq 2x/Q$ must in particular be $Px/Q$-smooth, this is all at most $\sum_{n \leq x} |a_n| \tau_{3}^{Px/Q}(n)$. Meanwhile, the Cauchy--Schwarz inequality and Non-concentration on Multiples condition yield $\sum_{2x/Q < h \leq Q} \frac{1}{h} \sum_{n_1 \leq x/h} |a_{n_1 h}| \leq \sqrt{\sum_{2x/Q < h \leq Q} \frac{1}{h^2}} \sqrt{\sum_{2x/Q < h \leq Q} ( \sum_{\substack{n \leq x, \\ h \mid n}} |a_{n}| )^2} \ll \sqrt{\frac{Q}{x K_{\text{conc}}[2x/Q]}} \sum_{n \leq x} |a_{n}|$.

Having replaced $a_{n_2 h}$ by $\frac{\sum_{n \leq x} a_n}{x}$ throughout \eqref{keyexprreplace}, we can use the Resembles the Integers I condition again to do the same for $a_{n_1 h}$. The resulting error terms are
\begin{eqnarray}
& \ll & \frac{(\sum_{n \leq x} |a_n|)^2}{x} \Biggl( \sum_{h \leq \frac{2x}{Q}} \sum_{0 \leq l \leq \frac{2x}{Qh}} \frac{\log(l+2)}{l+2} \sum_{\substack{m_1 \leq R, \\ P \; \text{smooth}}} \tau(m_1) \sum_{n_2 \leq x/h} \tau(\text{Sm}_{P}(n_2)) (\frac{Q}{K_{\text{int}}^{(1)} x m_1} + \frac{1}{K_{\text{int}}^{(2)}} ) + \nonumber \\
&& + \sum_{2x/Q < h \leq Q} \sum_{\substack{m_1 \leq R , \\ P \; \text{smooth}}} \sum_{n_2 \leq x/h} (\frac{1}{K_{\text{int}}^{(1)} hm_1} + \frac{1}{K_{\text{int}}^{(2)}} ) \Biggr) . \nonumber
\end{eqnarray}
Estimating similarly as before, and also using the bound $\sum_{n_2 \leq x/h} \tau(\text{Sm}_{P}(n_2)) \ll (x/h)\log P$, this is all readily seen to be
\begin{eqnarray}
& \ll & \frac{(\sum_{n \leq x} |a_n|)^2}{x} \Biggl( \log^{2}(2+\frac{x}{Q}) \log^{2}P \sum_{h \leq \frac{2x}{Q}} (\frac{Q \log P}{h K_{\text{int}}^{(1)}} + \frac{R x}{h K_{\text{int}}^{(2)}} ) + \sum_{\frac{2x}{Q} < h \leq Q} (\frac{x \log P}{h^2 K_{\text{int}}^{(1)}} + \frac{Rx}{h K_{\text{int}}^{(2)}} )  \Biggr) \nonumber \\
& \ll & (\sum_{n \leq x} |a_n|)^2 \Biggl( \frac{Q \log^{6}(P + \frac{x}{Q})}{x K_{\text{int}}^{(1)}} + \frac{R \log^{5}(P + \frac{x}{Q}) \log(Q^{2}/x)}{K_{\text{int}}^{(2)}} \Biggr) , \nonumber
\end{eqnarray}
which again is acceptable for Theorem \ref{likezthm}.

\vspace{12pt}
We have now shown that, up to acceptable error terms, the double sums in \eqref{keyexprreplace} are equal to $(\frac{\sum_{n \leq x} a_n}{x})^2$ multiplied by the corresponding double sums with all terms $a_n$ replaced by 1. Checking back to \eqref{keyexpression}, this means we have shown that
\begin{equation}\label{scalerel}
V(\mathcal{A},x,Q) - \frac{Q}{2} \sum_{n \leq x} a_n^2 = (\frac{\sum_{n \leq x} a_n}{x})^2 \left( V(\N,x,Q) - \frac{Q}{2} \sum_{n \leq x} 1 \right) + \text{error} ,
\end{equation}
where (with a small calculation, noting that in the error terms for $\mathcal{A} = \N$ we may take $K_{\text{prog}} \asymp Q \gg \sqrt{x}$ and $K_{\text{conc}}[1] \asymp 1$, as well as $K_{\text{int}}^{(3)} \asymp \log x$ and $K_{\text{int}}^{(2)} \asymp \sqrt{x}$ say) the overall ``error'' is still acceptably small for Theorem \ref{likezthm}.

To conclude, we now require some information about $V(\N,x,Q)$.

\begin{lem3}\label{nvarlem}
For all large $Q$ and $x$ we have
$$ V(\N,x,Q) = x^2 \int_{x/Q}^{2x/Q} \{v\} (1 - \{v\}) \frac{dv}{v^3} + O(Q\log^{2}Q + x) , $$
where $\{\cdot\}$ denotes the fractional part.
\end{lem3}

Note that, for completeness, we do not insist in Lemma \ref{nvarlem} that $Q \leq x$ (although that is the case in Theorem \ref{likezthm}).

\begin{proof}[Proof of Lemma \ref{nvarlem}]
It seems likely that this calculation has been performed in the literature before (perhaps many times, at least implicitly), but in the absence of a very obvious reference we provide a short proof.

First observe that for each $Q/2 < q \leq Q$, if we let $\N_{x,q} := \N \cap (q\lfloor \frac{x}{q} \rfloor,x]$ then we have
$$ \sum_{1 \leq a \leq q} \left|\N(x;q,a) - \text{Aver}(\N,x;q,(a,q)) \right|^{2} = \sum_{1 \leq a \leq q} \left|\N_{x,q}(x;q,a) - \text{Aver}(\N_{x,q},x;q,(a,q)) \right|^{2} ,$$
because the contributions to $\N(x;q,a)$ and $\text{Aver}(\N,x;q,(a,q))$ from those $n \leq q\lfloor \frac{x}{q} \rfloor$ come from complete sets of residue classes mod $q$, and so perfectly cancel. Then expanding the right hand side as in \eqref{stopensquare}, and noting that $\N_{x,q}(x;q,a)$ is always either 1 or 0, we find this is
$$ = \sum_{1 \leq a \leq q} \N_{x,q}(x;q,a)^{2} - \sum_{h \mid q} \phi(q/h) \text{Aver}(\N_{x,q},x;q,h)^{2} = \#\N_{x,q} - \sum_{h \mid q} \frac{1}{\phi(q/h)} ( \sum_{\substack{q \lfloor \frac{x}{q} \rfloor < n \leq x, \\ (n,q)=h}} 1 )^{2} . $$
Here we have $\#\N_{x,q} = x - q\lfloor \frac{x}{q} \rfloor + O(1) = q \{\frac{x}{q}\} + O(1)$, as well as the standard estimate $\sum_{\substack{q \lfloor \frac{x}{q} \rfloor < n \leq x, \\ (n,q)=h}} 1 = \sum_{\substack{q \lfloor \frac{x}{q} \rfloor/h < m \leq x/h, \\ (m,q/h)=1}} 1 = \frac{\phi(q/h)}{q/h} \sum_{q \lfloor \frac{x}{q} \rfloor/h < m \leq x/h} 1 + O(\tau(q/h)) = \phi(q/h) \{\frac{x}{q}\} + O(\tau(q/h))$ (see e.g. chapter 3.1 of Montgomery and Vaughan~\cite{mv} for the second equality here). Inserting these estimates, summing over $Q/2 < q \leq Q$, and calculating a little, we deduce that
$$ V(\N,x,Q) = \sum_{\frac{Q}{2} < q \leq Q} q \{\frac{x}{q}\} (1 - \{\frac{x}{q}\}) + O(Q + \sum_{\frac{Q}{2} < q \leq Q} \sum_{h \mid q} \tau(q/h)) = \sum_{\frac{Q}{2} < q \leq Q} q \{\frac{x}{q}\} (1 - \{\frac{x}{q}\}) + O(Q\log^{2}Q) . $$

Finally, note that if $x/(n+1) < t < x/n$ for some $n \in \N$ then the derivative of $t \mapsto t \{\frac{x}{t}\} (1 - \{\frac{x}{t}\}) = t (\frac{x}{t} - n) (n+1 - \frac{x}{t})$ is $O(x/t)$, and if $t > x$ then the derivative of $t \mapsto t \{\frac{x}{t}\} (1 - \{\frac{x}{t}\}) = x (1 - \frac{x}{t})$ is $x^{2}/t^2 \ll x/t$ as well. Thus we always have $q \{\frac{x}{q}\} (1 - \{\frac{x}{q}\}) = \int_{q}^{q+1} t \{\frac{x}{t}\} (1 - \{\frac{x}{t}\}) dt + O(\frac{x}{q})$, and summing over $Q/2 < q \leq Q$ yields
$$ V(\N,x,Q) = \sum_{\frac{Q}{2} < q \leq Q} \int_{q}^{q+1} t \{\frac{x}{t}\} (1 - \{\frac{x}{t}\}) dt + O(x + Q\log^{2}Q) = \int_{Q/2}^{Q} t \{\frac{x}{t}\} (1 - \{\frac{x}{t}\}) dt + O(x + Q\log^{2}Q) . $$
The Lemma now follows by substituting $v=x/t$ in the integral.
\end{proof}

Inserting the estimate from Lemma \ref{nvarlem} into \eqref{scalerel}, we obtain the claimed main terms in Theorem \ref{likezthm} along with an error term $O((\frac{\sum_{n \leq x} a_n}{x})^2 (x + Q\log^{2}Q)) = O(\frac{(\sum_{n \leq x} a_n)^2 \log^{2}x}{x})$. But this is smaller than the error term $O((\sum_{n \leq x} |a_n|)^2 \frac{\log x}{P(\log P) K_{\text{conc}}[1]})$ in the statement of Theorem \ref{likezthm}, noting that $P \leq x^{1/10}$ and that, by definition, $\frac{(\sum_{n \leq x} |a_n|)^2}{K_{\text{conc}}[1]} \geq (\sum_{n \leq x} |a_n|)^2$.
\qed

% SECTION 5 %%%%%%%%%%%%%%%%%%%%%%%%%%%%%
\section{Corollary \ref{smoothcor}: The smooth numbers example}
Before turning to Corollary \ref{smoothcor}, we recall a celebrated smooth numbers estimate of Hildebrand and Tenenbaum~\cite{ht}, that for all $2 \leq y \leq x$ we have
$$ \Psi(x,y) = \frac{x^{\alpha} \zeta(\alpha,y)}{\alpha \sqrt{2\pi(1+(\log x)/y)\log x \log y}} \left(1 + O\left(\frac{1}{\log(u+1)} + \frac{1}{\log y} \right) \right) . $$
Here $u = (\log x)/\log y$ and $\zeta(s,y) := \sum_{n \in \mathcal{S}(y)} \frac{1}{n^{s}} = \prod_{p \leq y}(1-p^{-s})^{-1}$, and $\alpha = \alpha(x,y) > 0$ is the ``saddle point'' corresponding to the $y$-smooth numbers less than $x$, which satisfies $\alpha(x,y) = 1 - \frac{\log(u\log(u+1))}{\log y} + O\left(\frac{1}{\log y}\right)$ provided $\log x < y \leq x$. In particular, when $\log^{C}x \leq y \leq x$ and $x$ is large we have $1-1/C + o(1) \leq \alpha(x,y) \leq 1 + o(1)$. On this range of $y$ we also have $\zeta(\alpha,y) = e^{O(u)}\log x$ (see e.g. Lemma 7.5 of Montgomery and Vaughan~\cite{mv}), and so $\Psi(x,y) = x e^{-u\log u - u\log\log(u+1) + O(u)}$.

\vspace{12pt}
We now begin by investigating the Progressions Condition. Similarly as discussed in the Introduction, in the smooth numbers case the left hand side there becomes
$$ \sum_{1 \leq r \leq 2x/Q} \sum_{\substack{h \mid r, \\ h \in \mathcal{S}(y)}} \sum_{\substack{a=1 , \\ (a,r) = h}}^{r} \Psi(x/h,y;r/h,a/h) \max_{z \leq x} \left|\Psi(z/h,y;r/h,a/h) - \frac{\Psi_{r/h}(z/h,y)}{\phi(r/h)} \right| . $$
Writing $r=Rh$ and proceeding a bit wastefully, this is all
$$ \leq \sum_{\substack{h \leq 2x/Q, \\ h \in \mathcal{S}(y)}} \sum_{R \leq 2x/Qh} \Psi_{R}(x/h,y) \max_{\substack{1 \leq a \leq R, \\ (a,R) = 1}} \max_{Z \leq x/h} \left|\Psi(Z,y;R,a) - \frac{\Psi_{R}(Z,y)}{\phi(R)} \right| . $$
Then Th\'eor\`eme 2.4(i) of La Bret\`eche and Tenenbaum~\cite{dlbten} implies that $\Psi_{R}(x/h,y) \leq \Psi(x/h,y) \ll \frac{\Psi(x,y)}{h^{\alpha}}$, where $\alpha = \alpha(x,y)$ is the saddle point. And using the Bombieri--Vinogradov type result\footnote{Theorem 1 of Harper~\cite{harpersmoothbv} is stated without the innermost maximum over $Z \leq x/h$, but that may be incorporated with a few modifications to the proof, and at the unimportant cost of a factor $\log^{9/2}x$ here rather than $\log^{7/2}x$ there. Indeed, we need only consider the maximum over $(x/h)^{0.99} \leq Z \leq x/h$, say, since for $Z \leq (x/h)^{0.99}$ the easy pointwise bound $\left|\Psi(Z,y;R,a) - \frac{\Psi_{R}(Z,y)}{\phi(R)} \right| \ll \frac{(x/h)^{0.99}}{\phi(R)}$ suffices. Then the zero-density part of the proof of Theorem 1 already uniformly bounds the maxima of the relevant character sums, after possibly adjusting the constants $c,C$. The large sieve part of that proof (in $\S\S 4.1-4.3$) may be adjusted to incorporate the maximum without changing the bounds, e.g. by breaking into dyadic ranges of $Z$ between $(x/h)^{0.99}$ and $x/h$; noting that the maximum over a dyadic range may be handled simply using the triangle inequality in the Perron separation of variables step of the argument; and then summing over the dyadic points. Finally, the part of the proof handling very large $y$ requires one extra Perron separation of variables (in the double sums over $m,n$ to which Vaughan's identity is applied) to incorporate the maximum, worsening the final bound by one logarithm.} for smooth numbers in Theorem 1 of Harper~\cite{harpersmoothbv}, provided that $y \geq \log^{C}x$ and $2x/Qh \leq \sqrt{\Psi(x/h,y)}$ we have
$$ \sum_{R \leq \frac{2x}{Qh}} \max_{\substack{1 \leq a \leq R, \\ (a,R) = 1}} \max_{Z \leq x/h} \left|\Psi(Z,y;R,a) - \frac{\Psi_{R}(Z,y)}{\phi(R)} \right| \ll_{A} \Psi(\frac{x}{h}, y)\left( \frac{e^{-\frac{cu}{\log^2(u+1)}}}{\log^{A}x} + y^{-c} \right) + \sqrt{\Psi(\frac{x}{h}, y)} \frac{x}{Qh} \log^{9/2}x . $$
Under our assumptions that $Q \geq x^{0.51}$ and $y \geq \log^{C}x$, we will indeed have $\sqrt{\Psi(x/h,y)} \geq \sqrt{(x/h)^{1-1/C+o(1)}} \geq 2x/Qh$ (provided $C$ is large enough), and the term $\Psi(\frac{x}{h}, y) ( \frac{e^{-\frac{cu}{\log^2(u+1)}}}{\log^{A}x} + y^{-c})$ in the Bombieri--Vinogradov bound is the dominant one. Thus the left hand side in the Progressions Condition is
$$ \ll_{A} \sum_{\substack{h \leq 2x/Q, \\ h \in \mathcal{S}(y)}} \Psi(\frac{x}{h}, y)^2 \Biggl( \frac{e^{-\frac{cu}{\log^2(u+1)}}}{\log^{A}x} + y^{-c} \Biggr) \ll \Psi(x, y)^2 \Biggl( \frac{e^{-\frac{cu}{\log^2(u+1)}}}{\log^{A}x} + y^{-c} \Biggr)  \sum_{h \leq 2x/Q} \frac{1}{h^{2\alpha}} , $$
and since the sum over $h$ is uniformly bounded we may take $\frac{1}{K_{\text{prog}}} \asymp_{A} ( \frac{e^{-\frac{cu}{\log^2(u+1)}}}{\log^{A}x} + y^{-c})$.

In the Non-concentration on Multiples condition, we can bound the left hand side by $\sum_{H \leq h \leq Q} \Psi(x/h,y)^2$. Again, Th\'eor\`eme 2.4(i) of La Bret\`eche and Tenenbaum~\cite{dlbten} implies this is $\ll \Psi(x,y)^2 \sum_{H \leq h \leq Q} \frac{1}{h^{2\alpha}}$, and so we may take $K_{\text{conc}}[H] \asymp H^{2\alpha - 1}$.

\subsection{The range $\log^{C}x \leq y \leq x^{1/\log\log x}$}
On this range we shall apply Theorem \ref{basicthm}, so we must also check the Hereditarily Sparse condition. The left hand side there is $\leq \sum_{n \in I \cap [1,x/h], n \in \mathcal{S}(y)} \tau(n)\ll \sum_{n \in I \cap [1,x/h], n \in \mathcal{S}(y)} \sum_{\substack{d|n, \\ d \leq \sqrt{x}}} 1 = \sum_{d \leq \sqrt{x}, d \in \mathcal{S}(y)} \sum_{md \in I\cap [1,x/h], m \in \mathcal{S}(y)} 1$. If the left endpoint of the interval $I$ is smaller than $Q/2$, we simply extend $I$ to the initial segment $[1,Q]$ and get a bound $\ll \sum_{d \leq \sqrt{x}, d \in \mathcal{S}(y)} \Psi(Q/d,y) \ll \Psi(x,y) (Q/x)^{\alpha} \sum_{d \leq \sqrt{x}, d \in \mathcal{S}(y)} \frac{1}{d^{\alpha}}$. Otherwise, the short interval estimate from Smooth Numbers Result 3 of Harper~\cite{harpersmoothrest} (say) is applicable, and implies $\sum_{md \in I\cap [1,x/h], m \in \mathcal{S}(y)} 1 \ll \Psi(x,y) (Q/xd)^{\alpha} \log x$. So we always get a bound $\ll \Psi(x,y) \log x (Q/x)^{\alpha} \sum_{d \leq \sqrt{x}, d \in \mathcal{S}(y)} \frac{1}{d^{\alpha}}$, and here $\sum_{d \leq \sqrt{x}, d \in \mathcal{S}(y)} \frac{1}{d^{\alpha}} \leq \zeta(\alpha,y) = e^{O(u)}\log x$.

Similarly, we find $\sum_{\substack{n \leq x, \\ n \in \mathcal{S}(y)}} \tau_3(n) = \sum_{\substack{d \leq x, \\ d \in \mathcal{S}(y)}} \tau(d) \Psi(x/d,y) \ll \Psi(x,y) \sum_{\substack{d \leq x, \\ d \in \mathcal{S}(y)}} \frac{\tau(d)}{d^{\alpha}}$ in the denominator on the right hand side of the Hereditarily Sparse condition. This is $\leq \Psi(x,y) \prod_{p \leq y} (1 - \frac{1}{p^{\alpha}})^{-2} \ll \Psi(x,y) e^{O(u)}\log^{2}x$. Thus we may finally take $\frac{1}{K_{\text{hered}}} \asymp (Q/x)^{\alpha} e^{O(u)} \log^{4}x$. When $y \leq x^{1/\log\log x}$ we have $u \geq \log\log x$, so the powers of $\log x$ in our bounds may be absorbed into the $e^{O(u)}$ terms.

So overall, Theorem \ref{basicthm} implies that on this range of $y$ we have
$$ V(\mathcal{S}(y),x,Q) = \frac{Q}{2} \Psi(x,y) + O_{A}\Biggl( \Psi(x,y)^2 \left(  \frac{e^{-\frac{cu}{\log^2(u+1)}}}{\log^{A}x} + y^{-c} + (\frac{Q}{x})^{2\alpha - 1} + (\frac{Q}{x})^{\alpha} e^{O(u)} \right) + \sum_{\substack{n \leq x, \\ n \in \mathcal{S}(y)}} \tau(n) \Biggr) . $$
The terms involving $Q/x$ are both $\ll Q\Psi(x,y) \frac{\Psi(x,y)}{x} (\frac{Q}{x})^{2(\alpha-1)} e^{O(u)}$. Recalling that $\alpha - 1 = -\frac{\log(u\log(u+1)) + O(1)}{\log y}$, that $x/Q \leq x^{0.49}$ (since we assume $Q \geq x^{0.51}$), and that $\frac{\Psi(x,y)}{x} = e^{-(1+o(1))u\log(u+1)}$, this bound is $\ll Q\Psi(x,y) e^{-(0.02+o(1))u\log(u+1)}$. Furthermore, the other ``main terms'' $- \frac{Q}{2} \frac{\Psi(x,y)^2}{x}, \Psi(x,y)^2 \int_{x/Q}^{2x/Q} \{v\} (1 - \{v\}) \frac{dv}{v^3}$ in Corollary \ref{smoothcor} are both $\ll Q\frac{\Psi(x,y)^2}{x} \ll Q\Psi(x,y) e^{-(1+o(1))u\log(u+1)}$. When $y \leq x^{1/\log\log x}$, and so $u \geq \log\log x$, these quantities may all be acceptably absorbed into the error terms in Corollary \ref{smoothcor}.

\subsection{The range $x^{1/\log\log x} < y \leq x$}
On this range we apply Theorem \ref{likezthm}. In the Resembles the Integers I condition, we claim that for any large $P \leq e^{\sqrt{\log y}}$ we may take $\frac{1}{K_{\text{int}}^{(1)}} \asymp \frac{\log(u+1) (\log P + (\log\log y)^{2.1}) \log P}{\log y} = \frac{u\log(u+1) (\log P + (\log\log y)^{2.1}) \log P}{\log x}$ and $K_{\text{int}}^{(2)} \asymp e^{(\log\log y)^{2.1}}$, say. Indeed, if $|I| \leq x e^{-2(\log\log y)^{2.1}}$ then the left hand side there is trivially $\ll x e^{-2(\log\log y)^{2.1}}$, which is certainly $\leq \frac{\Psi(x,y)}{e^{(\log\log y)^{2.1}}}$ when $y > x^{1/\log\log x}$. If $I \subseteq [1,x/h]$ is longer then necessarily $h \leq e^{2(\log\log y)^{2.1}} \leq y$, so on the left hand side in the Condition we will have
$$ \sum_{\substack{n \in I, \\ p|n \Rightarrow p > P}} \textbf{1}_{nh \in \mathcal{S}(y)} = \sum_{\substack{n \in I \cap \mathcal{S}(y), \\ p|n \Rightarrow p > P}} 1 = \sum_{n \in I \cap \mathcal{S}(y)} \sum_{d|n, d \in \mathcal{S}(P)} \mu(d) = \sum_{\substack{d \leq x, \\ d \in \mathcal{S}(P)}} \mu(d) \sum_{md \in I, m \in \mathcal{S}(y)} 1 . $$
Using Rankin's trick, the contribution here from all those $P^{(\log\log y)^{2.1}} < d \leq x$ is trivially $\ll \sum_{P^{(\log\log y)^{2.1}} < d \leq x, d \in \mathcal{S}(P)} \frac{x}{d} \ll \frac{x}{P^{2(\log\log y)^{2.1}/\log P}} \prod_{p \leq P} (1 - \frac{1}{p^{1-2/\log P}})^{-1} \ll \frac{x \log P}{e^{2(\log\log y)^{2.1}}}$, which again is $\leq \frac{\Psi(x,y)}{e^{(\log\log y)^{2.1}}}$ on this range of $y$. And with a bit of calculation, Theorems 1 and 3 and Lemma 1(v) of Hildebrand~\cite{hildebrand} imply that for $|I| > x e^{-2(\log\log y)^{2.1}}$ and $d \leq P^{(\log\log y)^{2.1}} \leq e^{(\log\log y)^{2.1} \sqrt{\log y}}$, we have $\sum_{md \in I, m \in \mathcal{S}(y)} 1 = (1 + O(\frac{(1 + \log(xd/|I|))\log(u+1)}{\log y})) \frac{\Psi(x,y)}{x} \sum_{md \in I} 1$, yielding that
$$ \sum_{\substack{n \in I, \\ p|n \Rightarrow p > P}} \textbf{1}_{nh \in \mathcal{S}(y)} = \frac{\Psi(x,y)}{x} \sum_{\substack{n \in I, \\ p|n \Rightarrow p > P}} 1 + O\Biggl(\frac{\log(u+1)}{\log y} \frac{\Psi(x,y) |I|}{x} \sum_{d \in \mathcal{S}(P)} \frac{\log(2xd/|I|)}{d} + \frac{\Psi(x,y)}{e^{(\log\log y)^{2.1}}} \Biggr) . $$
Since $\sum_{d \in \mathcal{S}(P)} \frac{\log(2xd/|I|)}{d} \ll (\log\log y)^{2.1} \sum_{d \in \mathcal{S}(P)} \frac{1}{d} + \sum_{d \in \mathcal{S}(P)} \frac{\log d}{d} \ll (\log\log y)^{2.1} \log P + \log^{2}P$, we see the claimed value of $\frac{1}{K_{\text{int}}^{(1)}}$ is indeed permissible.

In the Resembles the Integers II condition, we can simply bound the left hand side by $\sum_{n \in I} \tau(n)$, which is certainly $\ll \sum_{n \in I} \sum_{\substack{d|n, \\ d \leq \sqrt{x}}} 1 \ll |I|\log x + \sqrt{x}$. Thus we may take $K_{\text{int}}^{(3)} = e^{u\log u + u\log\log(u+1) + O(u)} \log x$, and $K_{\text{int}}^{(2)} = e^{(\log\log y)^{2.1}}$ as before.

\vspace{12pt}
Inserting these bounds into Theorem \ref{likezthm}, we may take $R = e^{(\log\log y)^{2.02}}$ and $P = e^{(\log\log y)^{1.01}}$, say. Then clearly the error terms there involving $\sum_{n \leq x} |a_n|^2 \tau(n), K_{\text{prog}}, K_{\text{conc}}[1]$ are acceptably small for Corollary \ref{smoothcor}, since $P$ is larger than any fixed power of $\log x$. We also (slightly crudely) have $\frac{RP}{K_{\text{int}}^{(2)}} + \frac{Q K_{\text{int}}^{(3)}}{x R^{1/\log P}} \ll \frac{e^{2(\log\log y)^{2.02}}}{K_{\text{int}}^{(2)}} + \frac{Q K_{\text{int}}^{(3)}}{x e^{(\log\log y)^{1.01}}} \ll \frac{1}{e^{(1/2)(\log\log y)^{1.01}}}$, and so the corresponding ``big Oh'' term in Theorem \ref{likezthm} is more than good enough. Note it was crucial here that $K_{\text{int}}^{(2)}$ is large compared with $RP$.

The term involving $K_{\text{conc}}[2x/Q]$ is $\ll \Psi(x,y)^2 \frac{\log P}{K_{\text{int}}^{(1)}} (\frac{Q}{x})^{\alpha} = Q \Psi(x,y) \frac{(\log\log y)^{1.01}}{K_{\text{int}}^{(1)}} \frac{\Psi(x,y)}{x} (\frac{Q}{x})^{\alpha - 1} \ll Q \Psi(x,y) \frac{(\log\log y)^{4.12} u\log(u+1)}{\log x} \frac{\Psi(x,y)}{x} (\frac{Q}{x})^{\alpha - 1}$. Since $x/Q \leq x^{0.49}$ and $\alpha - 1 = -\frac{\log(u\log(u+1)) + O(1)}{\log y}$, this is all $\ll Q\Psi(x,y) \frac{(\log\log y)^{4.12}}{\log x} e^{-0.5 u\log(u+1)}$, which again is acceptable. Notice this is rather delicate, because we do not get much saving from $\frac{1}{K_{\text{int}}^{(1)}}$.

Finally we must deal with the error term involving $\sum_{\substack{n \leq x, \\ n \in \mathcal{S}(y)}} \tau_3^{Px/Q}(n)$. If $x/Q \geq e^{(\log\log y)^{1.01}}$ then we can rely on our earlier bound $\sum_{\substack{n \leq x, \\ n \in \mathcal{S}(y)}} \tau_3^{Px/Q}(n) \leq \sum_{\substack{n \leq x, \\ n \in \mathcal{S}(y)}} \tau_3(n) \ll \Psi(x,y) e^{O(u)}\log^{2}x$, which is $\ll \Psi(x,y) \log^{O(1)}x$ on the present range of $y$. That error term is then $\ll \frac{Q}{x K_{\text{int}}^{(1)}} \Psi(x,y)^2 \log^{O(1)}x \ll \Psi(x,y)^2 e^{-(1/2)(\log\log y)^{1.01}}$, which suffices.

If $x/Q < e^{(\log\log y)^{1.01}}$, we have $\sum_{\substack{n \leq x, \\ n \in \mathcal{S}(y)}} \tau_3^{Px/Q}(n) = \sum_{\substack{d \leq x, \\ d \in \mathcal{S}(Px/Q)}} \tau(d) \Psi(x/d,y) \ll \Psi(x,y) \sum_{\substack{d \leq x, \\ d \in \mathcal{S}(Px/Q)}} \frac{\tau(d)}{d^{\alpha}} \leq \Psi(x,y) \prod_{p \leq Px/Q} (1 - \frac{1}{p^{\alpha}})^{-2}$. Since $\alpha = 1 -\frac{\log(u\log(u+1)) + O(1)}{\log y} \geq 1 - \frac{1}{\log(Px/Q)}$ (with very much room to spare), the product over primes here is $\ll \prod_{p \leq Px/Q} (1 - \frac{1}{p})^{-2} \ll \log^{2}(Px/Q)$, and our whole error term is $\ll \frac{Q}{x K_{\text{int}}^{(1)}} \Psi(x,y)^2 \log^{8}(Px/Q) \ll Q \Psi(x,y) \frac{(\log\log y)^{11.19} u\log(u+1)}{\log x} \frac{\Psi(x,y)}{x}$. As $\frac{\Psi(x,y)}{x} = e^{-u\log u - u\log\log(u+1) + O(u)}$, this bound too is acceptable for Corollary \ref{smoothcor}.
\qed

\vspace{12pt}
\noindent {\em Acknowledgements.} The author would like to thank Andrew Granville, for (long ago!) introducing him to the Barban--Davenport--Halberstam asymptotic problem for smooth numbers, and for his encouragement. He also thanks Ofir Gorodetsky for discussions and comments on a draft of this paper, and Bob Vaughan for his comments and for explaining his results on the smooth numbers problem.

\vspace{12pt}
\noindent {\em Rights.} For the purpose of open access, the author has applied a Creative Commons Attribution (CC-BY) licence to any Author Accepted Manuscript version arising from this submission.


\begin{thebibliography}{99}

\bibitem{barban} M. B. Barban. The ``large sieve'' method and its application to number theory. {\em Uspehi Mat. Nauk}, \textbf{21}, no. 1, pp 51-102. 1966 Translation in {\em Russian Math. Surveys}.

\bibitem{dlbfiormom} R. de la Bret\`{e}che, D. Fiorilli. Moments of moments of primes in arithmetic progressions. {\em Proc. Lond. Math. Soc. (3)}, \textbf{127}, no. 1, pp 165-220. 2023

\bibitem{dlbten} R. de la Bret\`{e}che, G. Tenenbaum. Propri\'{e}t\'{e}s statistiques des entiers friables. {\em The Ramanujan Journal}, \textbf{9}, pp 139-202. 2005

\bibitem{bruv} J. Br\"udern, R. C. Vaughan. A Montgomery--Hooley theorem for sums of two cubes. {\em Eur. J. Math.}, \textbf{7}, no. 4, pp 1616-1644. 2021

\bibitem{bruv2} J. Br\"udern, R. C. Vaughan. Sums of two unlike powers in arithmetic progressions. {\em Eur. J. Math.}, \textbf{8}, suppl. 1, pp S182-S213. 2022 

\bibitem{davhal} H. Davenport, H. Halberstam. Primes in arithmetic progressions. {\em Michigan Math. J.}, \textbf{13}, pp 485-489. 1966 {\em Corrigendum}: {\em Michigan Math. J.}, \textbf{15}, p 505. 1968

\bibitem{fiormartdis} D. Fiorilli, G. Martin. Disproving Hooley's conjecture. {\em J. Eur. Math. Soc. (JEMS)}, \textbf{25}, no. 12, pp 4791-4812. 2023

\bibitem{gallagher} P. X. Gallagher. The large sieve. {\em Mathematika}, \textbf{14}, pp 14-20. 1967

\bibitem{goldstonvaughan} D. A. Goldston, R. C. Vaughan. On the Montgomery--Hooley asymptotic formula. {\em Sieve methods, exponential sums, and their applications in number theory,} LMS Lecture Notes, vol. 237, Cambridge University Press, Cambridge, pp 117-142. 1997

\bibitem{granvillesmoothsurvey} A. Granville. Smooth numbers: computational number theory and beyond. {\em Algorithmic number theory: lattices, number fields, curves and cryptography}, Math. Sci. Res. Inst. Publ., vol. 44, Cambridge University Press, pp 267-323. 2008

\bibitem{harpersmoothbv} A. J. Harper. Bombieri--Vinogradov and Barban--Davenport--Halberstam type theorems for smooth numbers. Preprint available online at \url{http://arxiv.org/abs/1208.5992}

\bibitem{harpersmoothrest} A. J. Harper. Minor arcs, mean values, and restriction theory for exponential sums over smooth numbers. {\em Compos. Math.}, \textbf{152}, no. 6, pp 1121-1158. 2016

\bibitem{hildebrand} A. Hildebrand. On the number of positive integers $\leq x$ and free of prime factors $> y$. {\em J. Number Theory,} \textbf{22}, no. 3, pp 289-307. 1986

\bibitem{ht} A. Hildebrand, G. Tenenbaum. On integers free of large prime factors. {\em Trans. Amer. Math. Soc.,} \textbf{296}, no. 1, pp 265-290. 1986

\bibitem{hooley} C. Hooley. On the Barban-Davenport-Halberstam theorem. I. {\em J. Reine Angew. Math.,} \textbf{274/275}, pp 206-223. 1975

\bibitem{hooley10} C. Hooley. On the Barban-Davenport-Halberstam theorem: X. {\em Hardy-Ramanujan Journal}, \textbf{20}, pp 2-11. 1997

\bibitem{hooley9} C. Hooley. On the Barban-Davenport-Halberstam theorem: IX. {\em Acta Arithmetica}, \textbf{83}, pp 17-30. 1998

\bibitem{hooley14} C. Hooley. On the Barban-Davenport-Halberstam theorem: XIV. {\em Acta Arith.}, \textbf{101}, no. 3, pp 247-292. 2002

\bibitem{hooleysurv} C. Hooley. On Theorems of Barban-Davenport-Halberstam Type. In {\em Number Theory for the Millennium II}, edited by M. A. Bennett et al., published by A K Peters Ltd., Natick, MA. 2002

\bibitem{mast} D. Mastrostefano. A lower bound for the variance in arithmetic progressions of some multiplicative functions close to 1. Preprint available online at \url{https://arxiv.org/abs/2102.10589}

\bibitem{mont2} H. L. Montgomery. Primes in arithmetic progressions. {\em Michigan Math. J.}, \textbf{17}, pp 33-39. 1970

\bibitem{mv} H. L. Montgomery, R. C. Vaughan. {\em Multiplicative Number Theory I: Classical Theory.} First edition, published by Cambridge University Press. 2007

\bibitem{vaughangbdh1} R. C. Vaughan. On a variance associated with the distribution of general sequences in arithmetic progressions. I. {\em Phil. Trans. R. Soc. Lond. A}, \textbf{356}, pp 781-791. 1998

\bibitem{vaughangbdh} R. C. Vaughan. On a variance associated with the distribution of general sequences in arithmetic progressions. II. {\em Phil. Trans. R. Soc. Lond. A}, \textbf{356}, pp 793-809. 1998



\end{thebibliography}
\end{document}